\documentclass [reqno]{amsart}
\usepackage {latexsym,epsfig}
\usepackage [latin1]{inputenc}
\usepackage [all]{xy}

\newcommand {\OO}{{\mathcal O}}
\newcommand {\cc}{{\mathcal C}}
\newcommand {\LL}{{\mathcal L}}
\newcommand {\calP}{{\mathcal P}}
\newcommand {\MM}{{\mathcal M}}
\newcommand {\ZZ}{{\mathbb Z}}
\newcommand {\PP}{{\mathbb P}}
\newcommand {\QQ}{{\mathbb Q}}
\newcommand {\CC}{{\mathbb C}}

\newcommand {\vdim}{\mathop {\rm vdim}\,}
\newcommand {\codim}{\mathop {\rm codim}\,}
\newcommand {\PGL}{\mathop {\rm PGL}\,}
\newcommand {\Pic}{\mathop {\rm Pic}\,}
\newcommand {\isom}{\cong}
\newcommand {\xydiag}[1]{\[ \xymatrix @M=5pt {#1} \]}
\newcommand {\preprint}[2]{preprint \discretionary {#1/}{#2}{#1/#2}}

\newtheorem {theorem}{Theorem}[section]
\newtheorem {lemma}[theorem]{Lemma}
\newtheorem {proposition}[theorem]{Proposition}
\newtheorem {corollary}[theorem]{Corollary}
\theoremstyle {definition}
\newtheorem {definition}[theorem]{Definition}
\newtheorem {example}[theorem]{Example}
\theoremstyle {remark}
\newtheorem {remark}[theorem]{Remark}
\newtheorem {construction}[theorem]{Construction}
\newtheorem {convention}[theorem]{Convention}

\begin {document}


\title [Gromov-Witten invariants of very ample hypersurfaces]{Absolute and
  relative Gromov-Witten invariants of very ample hypersurfaces}
\author {Andreas Gathmann}
\address {Harvard University, Department of Mathematics, Science Center,
          1 Oxford Street, Cambridge, MA 02138, USA}
\email {andreas@math.harvard.edu}
\thanks {Funded by the DFG scholarship Ga 636/1--1.}

\begin {abstract}
  For any smooth complex projective variety $X$ and smooth very ample
  hypersurface $ Y \subset X $, we develop the technique of genus zero relative
  Gromov-Witten invariants of $Y$ in $X$ in algebro-geometric terms. We prove
  an equality of cycles in the Chow groups of the moduli spaces of relative
  stable maps that relates these relative invariants to the Gromov-Witten
  invariants of $X$ and $Y$. Given the Gromov-Witten invariants of $X$, we show
  that these relations are sufficient to compute all relative invariants, as
  well as all genus zero Gromov-Witten invariants of $Y$ whose homology and
  cohomology classes are induced by $X$.
\end {abstract}

\maketitle


Much work has been done recently on Gromov-Witten invariants related to
hypersurfaces. There are essentially two different problems that have been
studied. The first one is the question: how can one compute the Gromov-Witten
invariants of a hypersurface from those of the ambient variety
\cite {Be},\cite {G},\cite {K},\cite {LLY}? The second problem, mainly
studied from the point of view of symplectic geometry, is the theory of
relative Gromov-Witten invariants of a hypersurface \cite {IP1},\cite
{IP2},\cite {LR},\cite {R},\cite {V}. The goal of this paper is to show that
these two problems that have been studied completely independently so far are
in fact very closely related.

Let $X$ be a smooth complex projective variety and $ Y \subset X $ a smooth
very ample hypersurface. We start by giving a very short description of our
method to compute the genus zero Gromov-Witten invariants of $Y$ in terms of
those of $X$, skipping all technical details.

Fix $ n \ge 1 $ and $ \beta \in H_2 (X) $. For $ m \ge 0 $, we let $
\bar M_{(m)} $ (the official notation will be $ \bar M_{(m,0,\dots,0)}
(X,\beta) $) be a suitable compactification of the moduli space of all
irreducible stable maps $ (\PP^1,x_1,\dots,x_n,f) $ to $X$ such that $f$ has
multiplicity at least $m$ to $Y$ at the point $ x_1 $. Obviously, $ \bar
M_{(0)} $ should be just the ordinary moduli space of stable maps to $X$. On
the other hand, $ \bar M_{(Y \cdot \beta+1)} $ should correspond to the moduli
space of stable maps to $Y$, as all irreducible curves in $X$ having
multiplicity $ Y \cdot \beta +1 $ to $Y$ must actually lie inside $Y$.
Moreover, $ \bar M_{(m+1)} $ is a subspace of $ \bar M_{(m)} $ of (expected)
codimension one.

The strategy is now obvious: if we can describe the (virtual) divisor $ \bar
M_{(m+1)} $ in $ \bar M_{(m)} $ intersection-theoretically in terms of known
classes (and our main theorem \ref {main-thm} does precisely that), then we can
compute intersection products on $ \bar M_{(m+1)} $ if we can compute them on $
\bar M_{(m)} $. Iterating this procedure for $m$ from $0$ to $ Y \cdot \beta
$, this means that we can compute the Gromov-Witten invariants of $Y$ if we can
compute the Gromov-Witten invariants of $X$. In fact, we will show in a
forthcoming paper that this method reproves and generalizes the well-known
``mirror symmetry'' type formulas for Gromov-Witten invariants of certain
hypersurfaces \cite {Be},\cite {G},\cite {LLY}.

Let us make the step from multiplicity $m$ to $ m+1 $ a bit more precise. It is
easily seen that there is a section of a line bundle $ L_{(m)} $ on $ \bar
M_{(m)} $ whose zero locus describes exactly the condition that $f$ vanishes to
order at least $ m+1 $ along $Y$ at $ x_1 $. Hence one would naïvely expect
that $ \bar M_{(m+1)} $ is just the first Chern class of $ L_{(m)} $, which
turns out to be $ m \psi + ev^* Y $ (where $ \psi $ is the cotangent line class
and $ ev $ the evaluation map at the first marked point). However, this
intuition breaks down for those stable maps where $ x_1 $ lies on a component
that is completely mapped to $Y$ by $f$ (see the picture in construction
\ref {sigma-constr}), as $f$ actually has infinite multiplicity to $Y$ at
$ x_1 $ in this case. Thus we get correction terms from reducible curves of
that kind in our final equation. These correction terms are quite complicated,
but they can be recursively computed as they are made up of invariants of
smaller degree.

In this paper we will define more general spaces than the $ \bar M_{(m)} $
mentioned above. Namely, we allow the specification of multiplicities to $Y$
not only at the point $ x_1 $ but at all marked points. We call those moduli
spaces the spaces of relative stable maps, and equip them with virtual
fundamental classes. Intersection products on them are then called relative
Gromov-Witten invariants. Of course, they have the obvious (possibly virtual)
geometric interpretation as numbers of curves having given multiplicities to
$Y$ and satisfying some additional incidence conditions.

It should be said clearly that the specification of more than one multiplicity
is not necessary if one only wants to compute the Gromov-Witten invariants of
$Y$ from those of $X$. However, the general case fits nicely into the picture
and establishes the connection to the existing literature on relative
Gromov-Witten invariants, as these invariants have only been considered so far
in the case where the sum of the multiplicities is equal to $ Y \cdot \beta $
(i.e.\ where ``all intersection points with $Y$ are marked'').

The outline of the paper is as follows. In section \ref {sec-moduli} we define
the moduli spaces of relative stable maps and define their virtual fundamental
classes. The construction of the line bundles $ L_{(m)} $ and the moduli spaces
for the correction terms mentioned above is given in section \ref {sec-incr}.
At the end of this section we state our main theorem \ref {main-thm} that
describes how the moduli spaces of relative invariants change if one of the
multiplicities is increased by one. The proof of this theorem is done in two
steps. In the first step in section \ref {sec-proj} we look at the special case
where $ Y \subset X $ is a hyperplane in projective space. In this case no
virtual fundamental classes are needed, and the main theorem is established by
purely geometric analysis. The ideas for the main proofs of this section are
taken from \cite {V}. In the second step in section \ref {sec-general}, we
prove the general case by ``pulling back'' the result for hyperplanes in $
\PP^N $ along the morphism $ \bar M_n (X,\beta) \to \bar M_n (\PP^N,d) $
induced by the complete linear system $ |Y| $. Finally, in section \ref
{sec-enum} we prove that the main theorem can be used to compute the absolute
and relative Gromov-Witten invariants of $Y$ in terms of the Gromov-Witten
invariants of $X$. In a forthcoming paper, we will study the structure of these
computations and give some explicit examples.

A few remarks seem in order how this work is related to the existing
literature. The original ideas and motivation for our paper come from the work
of R. Vakil \cite {V}, who proved the main theorem under the following
restrictions: $ Y \subset X $ is a hyperplane in $ \PP^N $, the sum of the
prescribed multiplicities is equal to the degree of the curves, and one of the
multiplicities is raised from zero to one. It is interesting to note that he
used the main theorem in the opposite direction, namely to compute the
invariants of $X$ from those of $Y$. But the algorithm used there is very
specific to the case of a hyperplane in $ \PP^N $; it does not work for general
$ Y \subset X $.

All methods that have been known so far to compute Gromov-Witten invariants of
hypersurfaces $ Y \subset X $ need the existence of a torus action on $X$ and
use the techniques of equivariant cohomology and fixed point localization. In
the case where $Y$ is Calabi-Yau or Fano, the ``mirror symmetry'' results of
A. Givental \cite {G} and B. Lian et al.\ \cite {LLY} relate the Gromov-Witten
invariants of $Y$ to those of $X$ and express them in terms of certain
hypergeometric functions. Our methods are completely different; they do not
place any restrictions on the variety $X$ and do not require $Y$ to be
Calabi-Yau or Fano. In a forthcoming paper we will show that our equations
actually lead to the same hypergeometric functions as mentioned above.

Recently A. Bertram \cite {Be} has found another way to compute certain
Gromov-Witten invariants of Calabi-Yau and Fano hypersurfaces in projective
space. He also uses the torus action method, but does the actual computations
in a different way. It seems that his computations are closely related to ours,
but the exact relation to our methods is still unclear.

Relative Gromov-Witten invariants of any genus have been introduced in
symplectic geometry by A. Li and Y. Ruan \cite {LR} as well as E. Ionel and
T. Parker \cite {IP1},\cite {IP2}. They have been defined for any codimension
two symplectic submanifold $Y$ of a symplectic manifold $X$. The main
application in symplectic geometry is the splitting formula that expresses the
Gromov-Witten invariants of a symplectic sum $ X_1 \#_Y X_2 $ in terms of the
relative Gromov-Witten invariants of $Y$ in $ X_1 $ and $ X_2 $. E. Ionel has
informed me that \cite {IP2} together with the results announced in \cite {IP1}
can be used to prove a statement in the symplectic category that is analogous
to our main theorem.

The author would like to thank T. Graber, J. Harris, and R. Vakil for numerous
discussions. This work has been done at the Harvard University, to which the
author is grateful for hospitality.


\section {Moduli spaces of relative stable maps} \label {sec-moduli}

We begin with the description of the set-up and the definition of the moduli
spaces of relative stable maps. Let $X$ be a smooth complex projective variety
and $ Y \subset X $ a smooth very ample hypersurface. For notational
convenience, we denote by $ A^*(X) $ the ring of algebraic cohomology classes
of $X$ modulo torsion, and by $ H_2^+(X) $ the group of effective algebraic
homology classes of dimension two, modulo torsion.

Let $ \alpha = (\alpha_1,\dots,\alpha_n) $ be an $n$-tuple of non-negative
integers. As usual, for such an $n$-tuple we define $ |\alpha| := n $ and $
\sum \alpha := \sum_{i=1}^n \alpha_i $. If $ \alpha = (\alpha_1,\dots,\alpha_n)
$ and $ \alpha' = (\alpha'_1,\dots,\alpha'_m) $, we write $ \alpha \cup \alpha'
$ for $ (\alpha_1,\dots,\alpha_n,\alpha'_1,\dots,\alpha'_m) $. For $ 1 \le k
\le n $, we write $ \alpha \pm e_k $ for $ (\alpha_1,\dots,\alpha_k \pm 1,
\dots,\alpha_n) $.

Let $ n \ge 0 $ and let $ \beta \in H_2^+(X) $ be a non-zero homology class.
We denote by $ \bar M_n (X,\beta) := \bar M_{0,n} (X,\beta) $ the
Deligne-Mumford stack of $n$-pointed genus zero stable maps to $X$ of class $
\beta $ as defined in \cite {BM}.

The moduli space $ \bar M^Y_\alpha (X,\beta) $ that we want to construct should
be thought of as a compactification of the space of all irreducible stable maps
$ (\PP^1,x_1,\dots,x_n,f) $ to $X$ of class $ \beta $ that meet $Y$ in the
points $ x_i $ with multiplicity $ \alpha_i $ for all $i$. We define it first
as a subset of the set of geometric points of $ \bar M_n (X,\beta) $, but we
will see soon that it has the structure of a closed substack of $ \bar M_n
(X,\beta) $.

\begin {definition} \label {def-rel}
  With notations as above, we define $ \bar M^Y_\alpha (X,\beta) $ to be the
  locus in $ \bar M_n (X,\beta) $ of all stable maps $ (C,x_1,\dots,x_n,f) $
  such that
  \begin {enumerate}
  \item \label {x-inside}
    $ f(x_i) \in Y $ for all $i$ with $ \alpha_i > 0 $,
  \item \label {cycle-cond}
    $ f^* Y - \sum_i \alpha_i x_i \in A_0 (f^{-1}(Y)) $ is effective.
  \end {enumerate}
  If there is no risk of confusion we will write $ \bar M_\alpha (X,\beta)
  $ instead of $ \bar M^Y_\alpha (X,\beta) $.
\end {definition}

\begin {remark}
  Condition \ref {x-inside} is obviously necessary for \ref {cycle-cond} to
  make sense. The cycle class $ f^* Y \in A_0 (f^{-1}(Y)) $ is well-defined
  by \cite {F} chapter 6 as the intersection product $ Y \cdot C $ in $ Y
  \times_X C = f^{-1}(Y) $. Note that the Chow groups of a scheme are
  equal to the Chow groups of its underlying reduced scheme (see \cite {F}
  example 1.3.1 (a)), so we may replace $ f^{-1}(Y) $ by its underlying
  reduced scheme above. So, by abuse of notation, if we talk about connected
  (resp.\ irreducible) components of $ f^{-1}(Y) $ in the sequel we will always
  mean connected (resp.\ irreducible) components of the underlying reduced
  scheme of $ f^{-1}(Y) $.
\end {remark}

\begin {remark}
  For degree reasons, the space $ \bar M (X,\beta) $ is obviously empty if $
  \sum \alpha > Y \cdot \beta $, so we will tacitly assume from now on that $
  \sum \alpha \le Y \cdot \beta $.
\end {remark}

\begin {remark} \label {spell-out}
  The Chow group $ A_0 $ of a point as well as of (connected but not
  necessarily irreducible) genus zero curves is just $ \ZZ $, so condition
  \ref {cycle-cond} in definition \ref {def-rel} can be reformulated as
  follows: for any connected component $Z$ of $ f^{-1}(Y) $ we must have
  \begin {enumerate}
  \item \label {point-cond}
    if $Z$ is a point, it is either unmarked or a marked point $ x_i $ such
    that the multiplicity of $f$ at $ x_i $ along $Y$ is at least $ \alpha_i $,
  \item \label {curve-cond}
    if $Z$ is one-dimensional, let $ C^{(i)} $ for $ 1 \le i \le r $ be the
    irreducible components of $C$ not in $Z$ but intersecting $Z$, and let
    $ m^{(i)} $ be the multiplicity of $ f|_{C^{(i)}} $ at $ Z \cap C^{(i)} $
    along $Y$. Then we must have
      \[ Y \cdot f_* Z + \sum_{i=1}^r m^{(i)} \ge \sum_{x_i \in Z} \alpha_i. \]
  \end {enumerate}
\end {remark}

\begin {example}
  Let $ X = \PP^3 $, $ Y = H $ a plane, $ \beta = 5 \cdot [\mbox {line}] $, and
  $ \alpha = (1,2) $. In the following picture, the curve on the left is in $
  \bar M_{(1,2)} (X,\beta) $, whereas the one on the right is not (condition
  \ref {curve-cond} of remark \ref {spell-out} is violated for the line marked
  $Z$, as $ 1 + 1 \not\ge 2 + 1 $).
  \begin {center}
    \epsfig {file=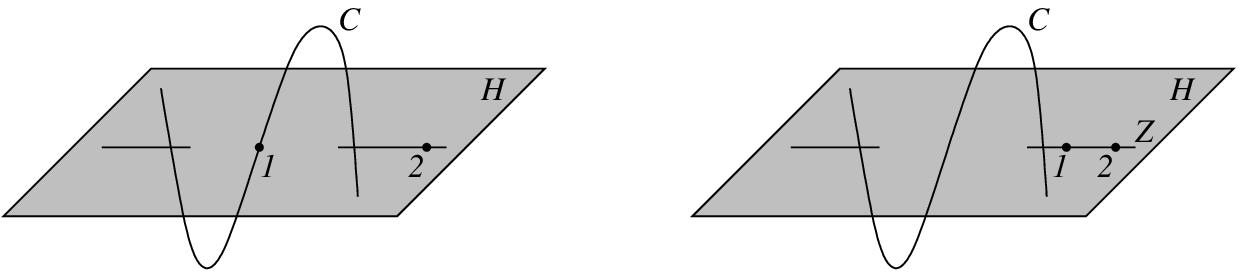}
  \end {center}
\end {example}

The first thing we will do is to study the space $ \bar M_\alpha (X,\beta) $ in
the special case where $ X = \PP^N $ and $ Y=H $ is a hyperplane. In this case,
we will write $ \bar M_\alpha (X,\beta) $ as $ \bar M_\alpha (\PP^N,d) $, where
$ d = H \cdot \beta $. The main result of this section is that the general
element of $ \bar M_\alpha (\PP^N,d) $ corresponds to an irreducible stable map
whose image is not contained in $H$, i.e.\ that the curves in $ \bar M_\alpha
(\PP^N,d) $ are exactly those that can be deformed to an \emph {irreducible}
curve that still satisfies the given multiplicity conditions and that is not
contained in $H$. (Here and in the following, by ``the curve $ \cc $ can be
deformed to a curve satisfying a property $P$'' we mean that there is a family
of stable maps such that the central fiber is $ \cc $ and the general fiber has
$P$.)

\begin {definition}
  We define $ M_\alpha (\PP^N,d) $ to be the subset of $ \bar M_\alpha
  (\PP^N,d) $ of all stable maps $ (C,x_1,\dots,x_n,f) $ with $ C \isom \PP^1 $
  and $ f(C) \not\subset H $.
\end {definition}

\begin {remark} \label {spell-out-sum}
  We will often consider first the easier case of the spaces $ \bar M_\alpha
  (\PP^N,d) $ with the additional condition that $ \sum \alpha = d $. (This is
  the situation that has been studied in \cite {V}.) In this case, condition
  \ref {cycle-cond} in definition \ref {def-rel} actually means that $ f^* H -
  \sum_i \alpha_i x_i = 0 \in A_0 (f^{-1}(H)) $. Correspondingly, the
  conditions in remark \ref {spell-out} read as follows: for any connected
  component $Z$ of $ f^{-1}(H) $ we must have
  \begin {enumerate}
  \item \label {point-cond-sum}
    if $Z$ is a point, it is a marked point $ x_i $ with $ \alpha_i $ being
    equal to the multiplicity of $f$ at $ x_i $ along $H$,
  \item \label {curve-cond-sum}
    if $Z$ is one-dimensional, let $ C^{(i)} $ for $ 1 \le i \le r $ be the
    irreducible components of $C$ not in $Z$ but intersecting $Z$, and let
    $ m^{(i)} $ be the multiplicity of $ f|_{C^{(i)}} $ at $ Z \cap C^{(i)} $
    along $H$. Then we must have
      \[ \deg f|_Z + \sum_{i=1}^r m^{(i)} = \sum_{x_i \in Z} \alpha_i. \]
  \end {enumerate}
\end {remark}

\begin {lemma} \label {loc-clos}
  The space $ M_\alpha (\PP^N,d) $ has the structure of an irreducible and
  locally closed substack of $ \bar M_n (\PP^N,d) $.
\end {lemma}

\begin {proof}
  The locus of irreducible stable maps $ (\PP^1,x_1,\dots,x_n,f) \in \bar M_n
  (\PP^N,d) $ such that $ f(\PP^1) \not\subset H $ can be written as $ M_n
  (\PP^N,d) \backslash \bar M_n (H,d) $, so it is open in $ \bar M_n (\PP^N,d)
  $. On the other hand, the condition that $f$ vanishes to order at least $
  \alpha_i $ along $H$ at $ x_i $ is closed, so $ M_\alpha (\PP^N,d) $ is the
  intersection of a closed subset with an open subset in $ \bar M_n (\PP^N,d)
  $. It is irreducible as there is a surjective rational map
    \[ \begin {array}{rcl}
         \CC^{2n} \times H^0 (\PP^1,\OO(d-\sum \alpha)) \times
                         H^0 (\PP^1,\OO(d))^N
         &\dashrightarrow& M_\alpha (\PP^n,d) \\
         (a_1,b_1,\dots,a_n,b_n,f_0,f_1,\dots,f_N)
         &\mapsto& (\PP^1,(a_1\!:\!b_1),\dots,(a_n\!:\!b_n),f)
       \end {array} \]
  where
    \[ f(z) = f(z_0:z_1) =
        (f_0(z) \cdot \prod_{i=1}^n (z_1 a_i - z_0 b_i)^{\alpha_i}:
          f_1(z):\cdots:f_N(z)) \]
  whose domain space is irreducible.
\end {proof}

\begin {lemma} \label {continuity}
  The closure of $ M_\alpha (\PP^N,d) $ in $ \bar M_n (\PP^N,d) $ is contained
  in $ \bar M_\alpha (\PP^N,d) $.
\end {lemma}

\begin {proof}
  This follows from the continuity of intersection products. To be more
  precise, let $ \cc $ be a point in the closure of $ M_\alpha (\PP^N,d) $. By
  lemma \ref {loc-clos} there is a family $ \phi: T \to \bar M_n (\PP^N,d) $ of
  stable maps over a smooth curve $T$ with a distinguished point $ 0 \in T $
  such that $ \phi(0) = \cc $ and $ \phi(t) \in M_\alpha (\PP^N,d) $ for $ t
  \neq 0 $. We have to prove that $ \phi(0) \in \bar M_\alpha (\PP^N,d) $.
  As it is obvious that $ \phi(0) $ satisfies condition \ref {x-inside} of
  definition \ref {def-rel}, it remains to show \ref {cycle-cond}.

  The family $ \phi $ is given by the data $ (C,x_1,\dots,x_n,f) $ where $ \pi:
  C \to T $ is a curve over $T$, the $ x_i: T \to C $ are sections of $ \pi $,
  and $ f: C \to \PP^N $ is a morphism. Set $ C_H = f^{-1}(H) $ and consider
  the 1-cycles $ f^* H $ and $ \sum_i \alpha_i x_i(T) $ in $ A_1(C_H) $. By
  assumption, the cycle $ \gamma := f^* H - \sum_i \alpha_i x_i(T) $ is
  effective (it might however have components over $ 0 \in T $ coming from
  $ f^* H $). Applying \cite {F} proposition 11.1 (b) to the cycles $ f^* H $
  and $ \gamma + \sum_i \alpha_i x_i(T) $ we see that the specialization
  of $ f^* H $ at $ t=0 $ is equal to the limit cycle of $ \gamma + \sum_i
  \alpha_i x_i(T) $ as $ t \to 0 $. As the limit cycle of $ \gamma $ for
  $ t \to 0 $ is effective, we have shown that $ \phi(0) $ satisfies \ref
  {cycle-cond}. This shows the lemma.
\end {proof}

\begin {definition}
  Let $ \cc = (C,x_1,\dots,x_n,f) \in \bar M_\alpha (\PP^N,d) $ be a stable
  map. An irreducible component $Z$ of $C$ is called an \emph {internal
  component} of $ \cc $ if $ f(C) \subset H $, and an \emph {external
  component} otherwise. A \emph {subcurve} of $ \cc $ is a stable map $ \cc' =
  (C',x_1',\dots,x_m',f') \in \bar M_{\alpha'} (\PP^N,d') $ constructed from $
  \cc $ as follows. Let $C'$ be any proper connected subcurve of $C$, and let $
  f'=f|_{C'} $. The marked points $ x_1',\dots,x_m' $ are the marked points $
  x_i $ contained in $C'$, together with all the intersection points of $C'$
  with the other irreducible components of $C$. We assign multiplicities $
  \alpha' = (\alpha_1',\dots, \alpha_m') $ to the points $ x_1',\dots,x_m' $ as
  follows: The points $ x_i $ on $C'$ will have their given multiplicity $
  \alpha_i $. The intersection points with other irreducible components of $C$
  will be assigned the multiplicity of $f'$ along $H$ at that point if the
  point lies on an external component of $C'$, and 0 otherwise. Let $d'$ be the
  degree of $f'$ on $C'$. The following picture shows an example of this
  construction, where the marked points are labeled with their multiplicities.
  \begin {center}
    \epsfig {file=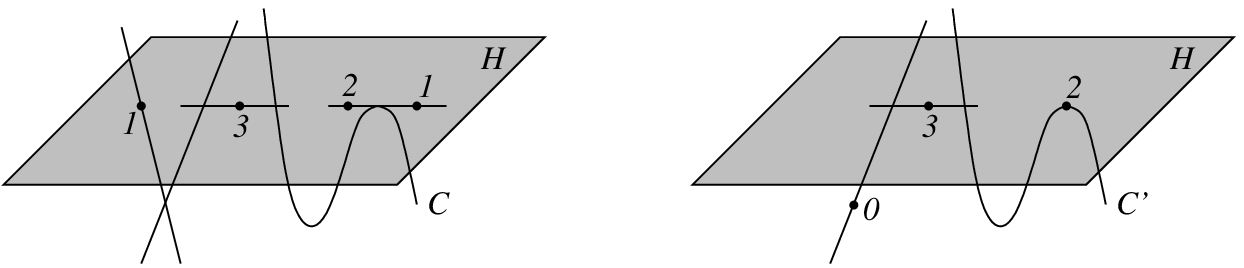}
  \end {center}
\end {definition}

\begin {lemma} \label {sum-ok}
  Let $ \cc \in \bar M_\alpha (\PP^N,d) $ be a stable map and assume that $
  \sum \alpha = d $. Let $ \cc' = (C',x_1',\dots,x_n',f') $ be a subcurve of $
  \cc $ with the following property: if $Z$ is an internal irreducible
  component of $C$ contained in $C'$, then any adjacent irreducible component
  of $Z$ in $C$ is also contained in $C'$. (For example, the subcurve in the
  picture above satisfies this property.) Then $ \sum \alpha' = d' $.
\end {lemma}

\begin {proof}
  The condition $ \sum \alpha = d $ means that $ f^* H - \sum \alpha_i x_i = 0
  \in A_0 (f^{-1}(H)) $. We claim that also $ {f'}^* H - \sum \alpha'_i x'_i =
  0 \in A_0 ({f'}^{-1}(H)) $, which then implies that $ \sum \alpha' = d' $. In
  fact, this can be checked on the connected components of $ {f'}^{-1}(H) $.
  Let $Z$ be a connected component of $ {f'}^{-1}(H) $. By assumption, there
  are only two possibilities:
  \begin {itemize}
  \item $C$ and $C'$ are locally isomorphic in a neighborhood of $Z$, i.e.\
    $Z$ is also a connected component of $ f^{-1}(H) $. Therefore, $ ({f'}^* H
    - \sum \alpha'_i x'_i)|_Z = 0 \in A_0 (Z) $.
  \item $Z$ is an intersection point of $C'$ with $ \overline {C \backslash C'}
    $ that lies on an external component of $C'$. Then, by definition of a
    subcurve, $Z$ is a marked point of $ \cc' $ with multiplicity equal to the
    multiplicity of $f'$ along $H$ at $Z$. In particular, we have again that
    $ ({f'}^* H - \sum \alpha'_i x'_i)|_Z = 0 \in A_0 (Z) $.
  \end {itemize}
  This proves the lemma.
\end {proof}

\begin {lemma} \label {def-basic}
  A stable map $ \cc = (C,x_1,\dots,x_n,f) \in \bar M_\alpha (\PP^N,d) $ can
  be deformed to an irreducible curve in $ \bar M_\alpha (\PP^N,d) $ if one
  of the following conditions is satisfied:
  \begin {enumerate}
  \item \label {int-int}
    $C$ has only internal components.
  \item \label {int-ext}
    $ \sum \alpha = d $, and $C$ consists exactly of one internal component $
    C^{(0)} $ and $r$ external components $ C^{(1)},\dots,C^{(r)} $
    intersecting $ C^{(0)} $ for some $ r \ge 0 $ (i.e.\ $C$ is a ``comb'',
    with the central component being internal and the teeth external, see the
    picture in construction \ref {sigma-constr}). Moreover, in this case
    $ \cc $ can even be deformed to an irreducible curve that is not contained
    in $H$ (which is then obvious unless $ r=0 $).
  \item \label {ext-ext}
    $ \sum \alpha = d $, and $C$ has exactly two irreducible components $
    C^{(1)} $ and $ C^{(2)} $, both being external.
  \end {enumerate}
\end {lemma}

\begin {proof}
  To show \ref {int-int}, note that by definition every curve with $ f(C)
  \subset H $ lies in $ \bar M_\alpha (\PP^N,d) $, so $ \bar M_n (H,d) \subset
  \bar M_\alpha (\PP^N,d) $. But it is well-known that the space of irreducible
  curves inside $ \bar M_n (H,d) $ is dense, so $ \cc $ can be deformed to an
  irreducible curve in $ \bar M_\alpha (\PP^N,d) $.

  \ref {int-ext} has been shown in \cite {V} theorem 6.1. (In fact, in the
  notations used in \cite {V}, our curve $ \cc $ is an element of a space $
  {\mathcal Y} $ with suitable decorations as introduced in \cite {V}
  definition 3.7.)

  Finally, in the case \ref {ext-ext} it is easy to construct an explicit
  deformation. Choose homogeneous coordinates $ z_0,\dots,z_N $ on $ \PP^N $
  such that $H$ is given by the equation $ z_0 = 0 $. The map $ f: C \to \PP^N
  $ is then given by sections $ s_0,\dots,s_N $ of a suitable line bundle $ \LL
  $ on $C$.  We may assume that the coordinates are chosen such that the $ s_i
  $ do not vanish at $ C^{(1)} \cap C^{(2)} $ (as for $ s_0 $ note that $ s_0
  (C^{(1)} \cap C^{(2)}) = 0 $ would mean that the intersection point lies on
  $H$, so it must be a marked point by remark \ref {spell-out-sum} \ref
  {point-cond-sum}, hence it must be non-singular, which is a contradiction).
  Let $ D_i = (s_i) $ be the associated divisors, in particular $ D_0 = \sum
  \alpha_i x_i $.

  Now let $W$ be the blow-up of $ \CC \times \PP^1 $ at the point $ (0,0) $,
  considered as a one-dimensional family of curves by the projection map $ \pi:
  W \to \CC $. We can identify the fiber $ \pi^{-1}(0) $ with $ C^{(1)} \cup
  C^{(2)} $. The points $ x_i \in \pi^{-1}(0) $ can be extended to sections $
  \tilde x_i $ of $ \pi $, giving rise to an extended divisor $ \tilde D_0 =
  \sum \alpha_i \tilde x_i $. In the same way one can find divisors $ \tilde
  D_i $ on $W$ such that $ \tilde D_i|_{\pi^{-1}(0)} = D_i $ for all $i$. As $
  \Pic W = \Pic C $, these divisors will be linearly equivalent and define a
  line bundle $ \tilde \LL $ on $W$ such that $ \tilde \LL|_{\pi^{-1}(0)} = \LL
  $. Moreover, after possibly restricting the base $ \CC $ to a smaller open
  neighborhood of $0$ we can assume that the $ \tilde D_i $ are base-point
  free. Finally, we can choose sections $ \tilde s_i $ of $ \tilde \LL $ such
  that $ (\tilde s_i) = \tilde D_i $ and $ \tilde s_i|{\pi^{-1}(0)} = s_i $.
  Then $ (W,\tilde x_0,\dots,\tilde x_n, (\tilde s_0:\cdots:\tilde s_N)) $ is a
  family of stable maps whose central fiber is $ \cc $ and whose general
  element is in $ M_\alpha (\PP^N,d) $.
\end {proof}

\begin {lemma} \label {def-glue}
  Let $ \cc = (C,x_1,\dots,x_n,f) \in \bar M_\alpha (\PP^N,d) $ be a reducible
  stable map and assume that $ \sum \alpha = d $. Then $ \cc $ can be deformed
  to a stable map in $ \bar M_\alpha (\PP^N,d) $ with fewer nodes.
\end {lemma}

\begin {proof}
  This is essentially obtained from lemma \ref {def-basic} by gluing. Pick a
  node $ P \in C $ and a subcurve $ \cc^{(0)} = (C^{(0)},x^{(0)}_1,\dots,
  x^{(0)}_{n^{(0)}},f^{(0)}) \in \bar M_{\alpha^{(0)}} (\PP^N,d^{(0)}) $ of $
  \cc $ as follows:
  \begin {enumerate}
  \item If $C$ has a node connecting two internal components of $C$, let $P$
    be this node and let $ C^{(0)} $ be the connected component of $ f^{-1}
    (H) $ containing $P$.
  \item Otherwise, if $C$ has a node connecting an internal component $Z$ to an
    external component of $C$, let $P$ be this node and let $ C^{(0)} $ be
    $Z$ together with all adjacent (necessarily external) components of $C$.
  \item Otherwise, let $P$ be any node of $C$ (necessarily connecting two
    external components of $C$) and let $ C^{(0)} $ be the two irreducible
    components of $C$ meeting at $P$.
  \end {enumerate}
  Let $ C^{(1)},\dots,C^{(r)} $ with $ r \ge 0 $ be the connected components
  of $ \overline {C \backslash C^{(0)}} $.

  In any case, we can deform $ \cc^{(0)} $ to an irreducible map in $ \bar
  M_{\alpha^{(0)}} (\PP^N,d^{(0)}) $ by lemma \ref {def-basic} (in the cases
  \ref {int-ext} and \ref {ext-ext} it follows from lemma \ref {sum-ok} that $
  \sum \alpha^{(0)} = d^{(0)} $). So let $ \phi: T \to \bar M_{\alpha^{(0)}}
  (\PP^N,d^{(0)}) $ be a deformation of $ \cc^{(0)} $ for some smooth pointed
  curve $ (T,0) $, i.e.\ $ \phi(0)=\cc^{(0)} $ and for all $ 0 \neq t \in T $
  the curve $ \phi(t) $ is irreducible. This deformation is given by a family $
  \pi: \tilde C \to T $ of curves, sections $ \tilde x_1,\dots,\tilde x_n $ of
  $ \pi $ and a map $ \tilde f: \tilde C \to \PP^N $. For all $ 1 \le i \le r
  $, the intersection point of $ C^{(0)} $ and $ C^{(i)} $ is one of the marked
  points of $ C^{(0)} $, hence corresponds to a marked point of $ \phi $, say $
  \tilde x_i $. Note that in all cases \ref {int-int} to \ref {ext-ext} above,
  the deformation $ \phi $ has the property that $ \tilde f (\tilde x_i(t)) \in
  H $ for all $ t \in T $ if this is true for $ t=0 $. In particular, there are
  $T$-valued projective automorphisms $ \psi_i: T \to \PGL (N) $ keeping $H$
  fixed such that $ \psi_i (t)(\tilde f(\tilde x_i(0))) = \tilde f (\tilde
  x_i(t)) $. The induced action of $ \PGL (N) $ on the moduli spaces $ \bar
  M_{\alpha^{(i)}} (\PP^N,d^{(i)}) $ makes $ \psi_i $ into a deformation of $
  \cc^{(i)} $ over $T$ such that for all $ t \in T $ the marked point
  corresponding to $ C^{(0)} \cap C^{(i)} $ is mapped to the same point in $
  \PP^N $ by the families $ \phi $ and $ \psi_i $. This means that the families
  $ \phi $ and $ \psi_i $ can actually be glued to give a deformation of the
  original curve $ \cc $. This deformation smoothes the node $P$.
\end {proof}

\begin {proposition} \label {closure}
  The closure of $ M_\alpha (\PP^N,d) $ in $ \bar M_n (\PP^N,d) $ is equal to $
  \bar M_\alpha (\PP^N,d) $. In particular, $ \bar M_\alpha (\PP^N,d) $ has the
  structure of an irreducible, proper, reduced substack of $ \bar M_n (\PP^N,d)
  $.
\end {proposition}

\begin {proof}
  ``$ \subset $'' has been shown in lemma \ref {continuity}, so it remains to
  show ``$ \supset $''. Let $ \cc \in \bar M_\alpha (\PP^N,d) $ be a stable
  map. Assume first that $ \sum \alpha = d $. Using lemma \ref {def-glue}
  inductively, we can deform $ \cc $ to an irreducible curve in $ \bar M_\alpha
  (\PP^N,d) $. If this irreducible curve does not lie inside $H$ then we are
  done, otherwise use the $ r=0 $ case of lemma \ref {def-basic} \ref
  {int-ext}.

  If $ k = d - \sum \alpha > 0 $, let $ \alpha' = \alpha \cup (1,\dots,1) $
  such that $ \sum \alpha' = d $. By adding marked points (and possibly
  introducing new contracted components) it is easy to find a stable map $
  \cc' \in \bar M_{\alpha'} $ that maps to $ \cc $ under the forgetful
  morphism $ \bar M_{n+k} (\PP^N,d) \to \bar M_n (\PP^N,d) $. By the above,
  $ \cc' $ can be deformed to an irreducible curve in $ M_{\alpha'} (\PP^N,d)
  $, which induces a deformation of $ \cc $ to an irreducible curve in $
  M_\alpha (\PP^N,d) $.

  Hence we finally have shown that $ \bar M_\alpha (\PP^N,d) $ is closed. So
  by giving it the reduced substack structure, we get a proper, reduced
  substack of $ \bar M_n (\PP^N,d) $ which is irreducible by lemma \ref
  {loc-clos}.
\end {proof}

\begin {lemma} \label {prop-pn}
  The moduli space $ \bar M_\alpha (\PP^N,d) $ has the following properties:
  \begin {enumerate}
  \item \label {fact-cov}
    If $ k = d - \sum \alpha > 0 $ and we let $ \alpha' = \alpha \cup
    (1,\dots,1) $ such that $ \sum \alpha' = d $, then there is a degree-$ k!
    $ generically finite cover $ \bar M_{\alpha'} (\PP^N,d) \to \bar M_\alpha
    (\PP^N,d) $, given by forgetting the last $k$ marked points and
    stabilizing.
  \item \label {univ-curve}
    $ \bar M_{\alpha \cup (0)} (\PP^N,d) $ is the universal curve over
    $ \bar M_\alpha (\PP^N,d) $. In particular, if $ \alpha = (0,\dots,0) $
    then $ \bar M_\alpha (\PP^N,d) = \bar M_{|\alpha|} (\PP^N,d) $.
  \item \label {exp-dim}
    The moduli space $ \bar M_\alpha (\PP^N,d) $ is purely of the expected
    dimension, which is $ \dim \bar M_{|\alpha|} (\PP^N,d) - \sum \alpha =
    d(N+1) + N - 3 + |\alpha| - \sum \alpha $.
  \end {enumerate}
\end {lemma}

\begin {proof}
  To show \ref {fact-cov}, note that from the parametrization of $ M_\alpha
  (\PP^N,d) $ given in the proof of lemma \ref {loc-clos} one can see that the
  general element of $ M_\alpha (\PP^N,d) $ corresponds to a stable map $
  (\PP^1,x_1,\dots,x_n,f) $ such that $ f^* H $ is equal to $ \sum_i \alpha_i
  x_i $ plus a union of $ k = d - \sum \alpha_i $ distinct unmarked points with
  multiplicity one. Obviously, the map $ \bar M_{\alpha'} (\PP^N,d) \to \bar
  M_\alpha (\PP^n,d) $ is finite over these elements, and it has degree $ k! $,
  corresponding to the choice of order of the $k$ added marked points.

  As in the proof of \ref {fact-cov}, the statement of \ref {univ-curve} is
  obvious on the dense open subset of $ \bar M_\alpha (\PP^N,d) $ described
  above, and it extends to the closures because of the flatness of the map $
  \bar M_{n+1} (\PP^N,d) \to \bar M_n (\PP^N,d) $.

  Finally, \ref {exp-dim} has been shown in \cite {V} proposition 5.7 if $
  \sum \alpha = d $. Otherwise use \ref {fact-cov} first. Alternatively, \ref
  {exp-dim} can be read off from the parametrization given in the proof of
  lemma \ref {loc-clos}.
\end {proof}

\begin {remark}
  The stack $ \bar M_\alpha (\PP^N,d) $ is in general singular, even in
  codimension one (see \cite {V} corollary 4.16). However, it is smooth at all
  points $ (\PP^1,x_1,\dots,x_n,f) \in M_\alpha (\PP^N,d) $. In fact, for these
  curves the obstruction space for deformations inside $ \bar M_\alpha
  (\PP^N,d) $ is $ H^1 (\PP^1,f^* T'_{\PP^N}) $, where $ f^* T'_{\PP^N} $ is
  the kernel of the composite map
    \[ f^* T_{\PP^N} \to f^* N_{H/\PP^N} \to (f^* N_{H/\PP^N})|_Z \]
  with $Z$ being the zero-dimensional subscheme of $ \PP^1 $ having length $
  \alpha_i $ at the point $ x_i $ for all $i$. But as both these maps are
  surjective on global sections (for the second one note that $ f^* N_{H/\PP^N}
  = \OO(d) $ and $ \sum \alpha \le d $), it follows that $ H^1 (\PP^1, f^*
  T'_{\PP^N}) = 0 $.

  However, we will not need any smoothness results in our paper.
\end {remark}

Now we return to the general case of the moduli space $ \bar M^Y_n (X,\beta) $
where $X$ is any smooth projective variety and $ Y \subset X $ a smooth
very ample hypersurface. One of the main problems is that these spaces will
in general not have the expected dimension. This means in particular that we
need virtual fundamental classes, which cannot be obtained using the techniques
above. To overcome this problem, we use the linear system $ |Y| $ to get a map
$ X \to \PP^N $, and consider the space $ \bar M^Y_\alpha (X,\beta) $ as the
``intersection'' of two problems we already know: (a) stable maps in $X$ and
(b) stable maps in $ \PP^N $ with given multiplicities to the hyperplane $ H
\subset \PP^N $ induced by $Y$.

We fix the following notation: let $ \varphi: X \to \PP^N $ be the morphism
determined by $ |Y| $, and let $ H \subset \PP^N $ the hyperplane such that $ Y
= \varphi^{-1}(H) $. As $ d := Y \cdot \beta > 0 $, the map $ \varphi $ induces
a morphism $ \phi: \bar M_n (X,\beta) \to \bar M_n (\PP^N,d) $ (see \cite
{BM}).

\begin {remark} \label {pullback-cond}
  Let $ \cc \in \bar M_n (X,\beta) $. As the conditions \ref {x-inside} and
  \ref {cycle-cond} of definition \ref {def-rel} pull back nicely, it is
  obvious that $ \cc \in \bar M^Y_\alpha (X,\beta) $ if and only if $ \phi
  (\cc) \in \bar M^H_\alpha (\PP^N,d) $.
\end {remark}

\begin {definition} \label {def-fund}
  By the previous remark, the space $ \bar M^Y_\alpha (X,\beta) $ has the
  structure of a proper closed substack of $ \bar M_n (X,\beta) $ by requiring
  the diagram of inclusions
  \xydiag {
    \bar M^Y_\alpha (X,\beta) \ar[r] \ar[d] &
      \bar M^H_\alpha (\PP^N,d) \ar[d] \\
    \bar M_n (X,\beta) \ar[r]^\phi &
      \bar M_n (\PP^N,d)
  }
  to be cartesian. We define the virtual fundamental class $ [\bar M^Y_\alpha
  (X,\beta)]^{virt} $ to be the one induced by the virtual fundamental class of
  $ \bar M_n (X,\beta) $ (see e.g.\ \cite {B},\cite {BF}) and the usual
  fundamental class of $ \bar M^H_\alpha (\PP^N,d) $, in the sense of the
  following remark.
\end {definition}

\begin {remark} \label {ind-fund}
  Let $ M_1 $ and $ M_2 $ be Deligne-Mumford stacks over a smooth
  Deligne-Mumford stack $S$. Let $ M = M_1 \times_S M_2 $, so that we have
  a cartesian diagram
  \xydiag {
    M \ar[r] \ar[d] & M_1 \times M_2 \ar[d] \\
    S \ar[r]^-{\Delta} & S \times S.
  }
  Assume that we are given classes $ \gamma_1 \in A_*(M_1) $ and $ \gamma_2 \in
  A_*(M_2) $ (usually thought of as virtual fundamental classes in this paper).
  Then the class $ \Delta^! (\gamma_1 \otimes \gamma_2) $ in $M$ will be called
  induced by $ \gamma_1 $ and $ \gamma_2 $. If the maps $ M_1 \to S $ and $ M_2
  \to S $ are inclusions, this is actually the usual refined intersection
  product of $ \gamma_1 $ and $ \gamma_2 $. This is the case in the above
  definition, but we mentioned the general case here as we will need it later
  on.
\end {remark}

By lemma \ref {prop-pn} \ref {exp-dim}, the virtual fundamental class of $
\bar M^Y_\alpha (X,\beta) $ defined above has dimension $ \dim \bar M_n (X,
\beta) - \sum \alpha $, which is the expected dimension of $ \bar M^Y_\alpha
(X,\beta) $. If $X$ is a projective space and $ Y \subset X $ a hyperplane, it
is obvious by definition that the virtual fundamental class of $ \bar
M^Y_\alpha (X,\beta) $ is equal to the usual one.


\section {Increasing the multiplicities} \label {sec-incr}

By construction, $ \bar M_{\alpha+e_k} (X,\beta) $ is a closed substack of $
\bar M_\alpha (X,\beta) $ of expected codimension one. The main goal of this
paper is to compute $ [\bar M_{\alpha+e_k} (X,\beta)]^{virt} $ as a cycle in
the Chow group of $ \bar M_\alpha (X,\beta) $. We start with the following
naïve approach describing the transition from multiplicity $ \alpha_k $ to $
\alpha_k + 1 $ at the point $ x_k $.

\begin {construction} \label {sigma-constr}
  Consider a moduli space $ M = \bar M_n (X,\beta) $ and let $ \cc \to M $ be
  the universal curve, with evaluation map $ ev: \cc \to X $. Fix $k$ with $ 1
  \le k \le n $ and let $ s_k: M \to \cc $ denote the section corresponding to
  the marked point $ x_k $. Let $ y \in H^0(\OO_X(Y)) $ be the equation of $Y$.
  Choose an integer $ m \ge 0 $. We pull $y$ back to $ \cc $ by $ ev $, take
  the $m$-jet relative to $M$ of it and pull this back to $M$ by $ s_k $ to
  get a section
    \[ \sigma_k^m := s_k^* d^m_{\cc/M} ev^* y \in
         H^0 (M,s_k^* \calP^m_{\cc/M} (ev^* \OO_X(Y)) ), \]
  where $ \calP^m_{\cc/M} (ev^* \OO_X(Y)) $ denotes relative principal parts of
  order $m$ (or $m$-jets) of the line bundle $ ev^* \OO_X(Y) $, and $
  d^m_{\cc/M} $ is the derivative up to order $m$ (see \cite {EGA4} 16.3,
  16.7.2.1 for precise definitions). Geometrically, $ \sigma_k^m $ vanishes
  precisely on the stable maps that have multiplicity at least $ m+1 $ to $Y$
  at the point $ x_k $. By \cite {EGA4} 16.10.1, 16.7.3 there is an exact
  sequence
    \[ 0 \to L_k^{\otimes m} \otimes ev_k^* \OO_X(Y)
         \to s_k^* \calP^m_{\cc/M} (ev^* \OO_X(Y))
         \to s_k^* \calP^{m-1}_{\cc/M} (ev^* \OO_X(Y))
         \to 0 \]
  where we set $ \calP^{-1}_{\cc/M}(ev^* \OO_X(Y)) = 0 $, and where $ L_k =
  s_k^* \omega_{\cc/M} $ is the $k$-th cotangent line, i.e.\ the line bundle
  on $M$ whose fiber at a point $ (C,x_1,\dots,x_n,f) $ is $ T^\vee_{C,x_k} $.
  Note that the last map in this sequence sends $ \sigma_k^m $ to $
  \sigma_k^{m-1} $ for $ m>0 $. Now restrict these bundles and sections to $
  \bar M_\alpha (X,\beta) $. As all stable maps in $ \bar M_\alpha (X,\beta) $
  have multiplicity (at least) $ \alpha_k $ at $ x_k $, the restriction of $
  \sigma_k^{\alpha_k} $ to $ \bar M_\alpha (X,\beta) $ defines a section
    \[ \sigma_k := \sigma_k^{\alpha_k}|_{\bar M_\alpha (X,\beta)} \in
         H^0 (L_k^{\otimes \alpha_k} \otimes ev_k^* \OO_X(Y))
       = H^0 (\OO (\alpha_k \psi_k + ev_k^* Y)) \]
  on $ \bar M_\alpha (X,\beta) $, where $ \psi_k = c_1 (L_k) $.
\end {construction}

The vanishing of this section describes exactly the condition that a stable map
in $ \bar M_\alpha (X,\beta) $ vanishes up to order $ \alpha_k+1 $ at $ x_k $.
Hence naïvely one would expect that $ \bar M_{\alpha+e_k} (X,\beta) $ is
described inside $ \bar M_\alpha (X,\beta) $ by the vanishing of this section,
and that $ [\bar M_{\alpha+e_k} (X,\beta)]^{virt} $ is given by
\begin {equation} \label {eqn-virt}
  (\alpha_k \, \psi_k + ev_k^* Y) \cdot [\bar M_\alpha (X,\beta)]^{virt}.
\end {equation}
This is not true, however, because of the presence of stable maps with the
property that the component on which $ x_k $ lies is mapped entirely into $Y$.
Of course, the section $ \sigma_k $ vanishes on those stable maps, but they are
in general not in $ \bar M_{\alpha+e_k} (X,\beta) $. Hence, these stable maps
will also contribute to the expression (\ref {eqn-virt}). We will now introduce
the moduli spaces of the stable maps occurring in these correction terms.
Informally speaking, generic stable maps in these correction terms have $ r+1 $
irreducible components $ C^{(0)},\dots,C^{(r)} $ for some $ r \ge 0 $, where
$ C^{(0)} $ (called the internal component) is mapped into $Y$, and all $
C^{(i)} $ for $ i>0 $ (called the external components) intersect $ C^{(0)} $
and have a prescribed multiplicity $ m^{(i)} $ to $Y$ at this intersection
point (see the picture below, where $ m^{(1)}=1 $ and $ m^{(2)}=2 $). The
point $ x_k $ has to lie on $ C^{(0)} $. The initial multiplicity conditions
$ \alpha $ as well as the homology class $ \beta $ get distributed in all
possible ways to the components $ C^{(i)} $.
\begin {center}
  \epsfig {file=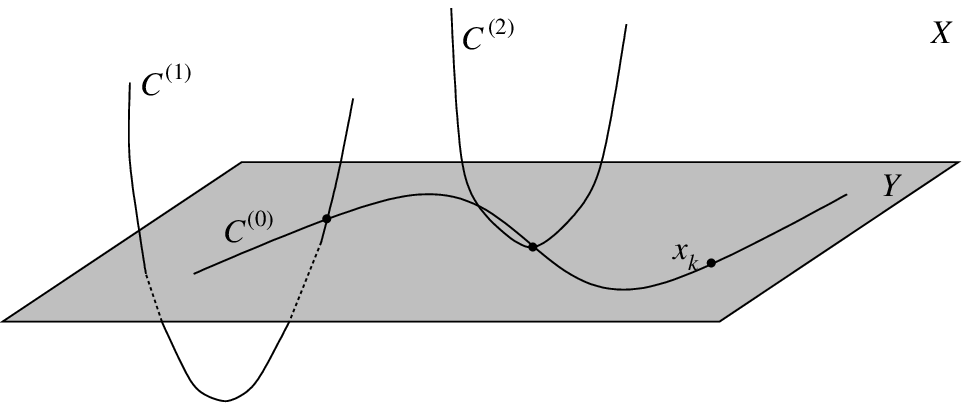}
\end {center}
We now describe this more formally.

\begin {definition} \label {def-dc}
  Consider a moduli space $ \bar M_\alpha (X,\beta) $ and $ 1 \le k \le n $ as
  above. Let $r$ be a non-negative integer. Choose a partition $
  A=(\alpha^{(0)},\dots,\alpha^{(r)}) $ of $ \alpha $ such that $ \alpha_k \in
  \alpha^{(0)} $. Let $ B=(\beta^{(0)},\dots,\beta^{(r)}) $ be an $(r+1)$-tuple
  of homology classes with $ \beta^{(0)} \in H_2^+(Y) $ and $ \beta^{(i)} \in
  H_2^+(X) \backslash \{0\} $ for $ i>0 $ such that $ i_* \beta^{(0)} + \beta^{
  (1)} + \cdots + \beta^{(r)} = \beta $, where $ i: Y \to X $ is the
  inclusion. Finally, choose an $r$-tuple $ M=(m^{(1)},\dots,m^{(r)}) $ of
  positive integers. With these notations, we define the moduli space $ D_k
  (X,A,B,M) $ to be the fiber product
    \[ D_k(X,A,B,M) := \bar M_{|\alpha^{(0)}|+r} (Y,\beta^{(0)}) \times_{Y^r}
         \prod_{i=1}^r \bar M_{\alpha^{(i)} \cup (m^{(i)})} (X,\beta^{(i)}) \]
  where the map from the first factor to $ Y^r $ is the evaluation at the last
  $r$ marked points, and the map from the second factor to $ Y^r $ is the
  evaluation at the last marked point of each of its factors. We define the
  virtual fundamental class of $ D_k(X,A,B,M) $ to be $ \frac {m^{(1)} \cdots
  m^{(r)}}{r!} $ times the class induced by the virtual fundamental classes of
  its factors, in the sense of remark \ref {ind-fund}. The reason for the
  unusual multiplicity will become clear in the proof of proposition \ref
  {main-p1-sum}.
\end {definition}

\begin {definition} \label {def-d}
  With notations as in the previous definition, let $ D_{\alpha,k} (X,\beta) $
  be the disjoint union of the $ D_k(X,A,B,M) $ for all possible $A$, $B$, and
  $M$ satisfying
  \begin {equation} \label {eqn-dim}
    Y \cdot i_* \beta^{(0)} + \sum_i m^{(i)} = \sum \alpha^{(0)}
  \end {equation}
  (the reason for this condition will become clear in the following lemma).
  The virtual fundamental class of $ D_{\alpha,k} (X,\beta) $ is defined to be
  the sum of the virtual fundamental classes of its components $ D_k(X,A,B,M)
  $.
\end {definition}

\begin {lemma} \label {irr-d-pn}
  In the case where $ X=\PP^N $ and $ Y=H $ is a hyperplane, the moduli spaces
  $ D_k(\PP^N,A,B,M) $ satisfying equation (\ref {eqn-dim}) of definition \ref
  {def-d} are proper irreducible substacks of $ \bar M_\alpha (\PP^N,d) $ of
  codimension one.
\end {lemma}

\begin {proof}
  Considering the definition of the space $ D_k(X,A,B,M) $, the fact that it
  is irreducible follows from the following three observations:
  \begin {enumerate}
  \item $ \bar M_{|\alpha^{(0)}|+r} (H,d^{(0)}) $ is irreducible,
  \item \label {flat-surj}
    the evaluation maps $ \bar M_{\alpha^{(i)} \cup (m^{(i)})} (\PP^N,d^{(i)})
    \to H $ at the last marked point are flat and surjective (this follows
    from the action of the group of automorphisms of $ \PP^N $ keeping $H$
    fixed on the space $ \bar M_{\alpha^{(i)} \cup (m^{(i)})} (\PP^N,d^{(i)})
    $),
  \item the fibers of the maps in \ref {flat-surj} are irreducible (by the
    Bertini theorem, as the spaces $ \bar M_{\alpha^{(i)} \cup (m^{(i)})}
    (\PP^N,d^{(i)}) $ itself are irreducible by proposition \ref {closure}).
  \end {enumerate}
  Moreover, these arguments show that the dimension of $ D_k(\PP^N,A,B,M) $
  is equal to
    \[ \dim \bar M_{|\alpha^{(0)}|+r} (H,d^{(0)})
       + \sum_{i=1}^r \dim \bar M_{\alpha^{(i)} \cup (m^{(i)})} (\PP^N,d^{(i)})
       - r \cdot (N-1). \]
  By a quick computation using lemma \ref {prop-pn} \ref {exp-dim} this
  is equal to
    \[ \dim \bar M_\alpha (\PP^N,d)
       + \sum \alpha^{(0)} - d^{(0)} - \sum_i m^{(i)} - 1, \]
  so the dimension statement follows from equation (\ref {eqn-dim}) of
  definition \ref {def-d}.

  The stack $ D_k(\PP^N,A,B,M) $ is visibly a closed substack of
    \[ \bar M_{|\alpha^{(0)}|+r} (\PP^N,d^{(0)}) \times_{(\PP^N)^r}
         \prod_{i=1}^r \bar M_{|\alpha^{(i)}|+1} (\PP^N,d^{(i)}), \]
  which in turn is a closed substack of $ \bar M_n (\PP^N,d) $ by \cite {BM}
  chapter 7 property III. To prove that it is contained in $ \bar M_\alpha
  (\PP^N,d) $ it suffices to show that a general element $ \cc = (C,x_1,\dots,
  x_n,f) \in D_k(\PP^N,A,B,M) $ satisfies the conditions of remark \ref
  {spell-out}. As $ \cc $ is general, we have $ C = C^{(0)} \cup \cdots \cup
  C^{(r)} $ where $ C^{(0)} \in M_{r+|\alpha^{(0)}|} (H,d^{(0)}) $ and
  $ C^{(i)} \in M_{\alpha^{(i)} \cup (m^{(i)})} (\PP^N,d^{(i)}) $. The
  condition of remark \ref {spell-out} is obvious for all connected components
  of $ f^{-1}(H) $ besides $ C^{(0)} $. As for $ C^{(0)} $, the condition is
  exactly the ``$ \ge $'' part of equation (\ref {eqn-dim}) of definition \ref
  {def-d}.
\end {proof}

\begin {remark}
  We will see in proposition \ref {pullback-d} that even for general $X$, the
  moduli spaces $ D_k(X,A,B,M) $ satisfying equation (\ref {eqn-dim}) of
  definition \ref {def-d} are proper substacks of $ \bar M_\alpha (X,\beta) $
  of expected codimension one. Thus we can view the virtual fundamental class
  of the $ D_k(X,A,B,M) $ as well as of $ D_{\alpha,k} (X,\beta) $ as cycles
  in the Chow group of $ \bar M_\alpha (X,\beta) $ whose dimension is equal to
  the expected dimension of $ \bar M_\alpha (X,\beta) $ minus one.
\end {remark}

We can now state the main theorem of this paper.

\begin {theorem} \label {main-thm}
  With notations as above, we have
    \[ (\alpha_k \, \psi_k + ev_k^* Y) \cdot [\bar M_\alpha (X,\beta)]^{virt}
         = [\bar M_{\alpha+e_k} (X,\beta)]^{virt}
           + [D_{\alpha,k} (X,\beta)]^{virt} \]
  in the Chow group of $ \bar M_\alpha (X,\beta) $, for all $ 1 \le k \le n $.
\end {theorem}

The proof will be given at the end of section \ref {sec-general}.


\section {Proof of the main theorem for hyperplanes in $ \PP^N $}
  \label {sec-proj}

In this section we will prove the main theorem \ref {main-thm} in the case
where $ X=\PP^N $ and $ Y=H $ is a hyperplane. Most of the proofs are
generalized versions from those in \cite {V}, where the generalizations are
quite straightforward. Recall that in construction \ref {sigma-constr} we
defined a section $ \sigma_k $ of a suitable line bundle on $ \bar M_\alpha
(\PP^N,d) $ such that the zero locus of $ \sigma_k $ has class $ \alpha_k \,
\psi_k + ev_k^* H $ and describes exactly those stable maps $ (C,x_1,\dots,
x_n,f) $ where $f$ vanishes to order at least $ \alpha_k+1 $ along $H$ at $
x_k $. For simplicity, we will restrict ourselves first to the case $ \sum
\alpha = d $ (note that the term $ [\bar M_{\alpha+e_k} (\PP^N,d)]^{virt} $
in the main theorem is then absent for degree reasons). We begin by proving a
set-theoretic version of the main theorem.

\begin {lemma} \label {sigma-set}
  Assume that $ \sum \alpha = d $. Then the zero locus of the section $
  \sigma_k $ on $ \bar M_\alpha (\PP^N,d) $ is equal to $ D_{\alpha,k}
  (\PP^N,d) $.
\end {lemma}

\begin {proof}
  By construction, it is obvious that $ \sigma_k $ vanishes on $ D_{\alpha,k}
  (\PP^N,d) $, so let us prove the converse. Let $ \cc = (C,x_1,\dots,x_n,f)
  \in \bar M_\alpha (\PP^N,d) $ be a stable map with $ \sigma_k (\cc) = 0 $.

  Assume first that $ x_n $ is an isolated point of $ f^{-1}(H) $. As $f$
  vanishes to order at least $ \alpha_k+1 $ along $H$ at $ x_k $, this is a
  contradiction to remark \ref {spell-out-sum} \ref {point-cond-sum}.

  So $ x_n $ is not an isolated point of $ f^{-1}(H) $. Let $ C^{(0)} $ be
  the connected component of $ f^{-1}(H) $ containing $ x_k $, and let $
  C^{(1)},\dots,C^{(r)} $ be the connected components of $ \overline {C
  \backslash C^{(0)}} $. Let $ m^{(i)} $ be the multiplicity of $ f|_{C^{(i)}}
  $ at $ C^{(0)} \cap C^{(i)} $ along $H$, let $ d^{(i)} $ be the degree of
  $f$ on $ C^{(i)} $, and let $ \alpha^{(i)} $ be the collection of the
  multiplicities $ \alpha_j $ such that $ x_j \in C^{(i)} $. Then it is obvious
  that $ \cc \in D_k (\PP^N,A,B,M) $ with $A$, $B$, $M$ as in definition \ref
  {def-dc}. Moreover, equation (\ref {eqn-dim}) of definition \ref {def-d} is
  satisfied by remark \ref {spell-out-sum} \ref {curve-cond-sum} applied to
  $ C^{(0)} $, hence it follows that $ \cc \in D_{\alpha,k} (\PP^N,d) $.
\end {proof}

\begin {remark} \label {factorial}
  As the spaces $ D_k(\PP^N,A,B,M) $ are irreducible and of codimension one by
  lemma \ref {irr-d-pn}, lemma \ref {sigma-set} tells us that in the case $
  \sum \alpha = d $ we must have
    \[ (\alpha_k \, \psi_k + ev_k^* H) \cdot [\bar M_\alpha (\PP^N,d)]
         = \sum \lambda_{A,B,M} \; [D_k(\PP^N,A,B,M)]^{virt} \]
  for some $ \lambda_{A,B,M} $, where the sum is taken over all $ A,B,M $
  for which $ D_k(\PP^N,A,B,M) $ occurs in $ D_{\alpha,k} (\PP^N,d) $. Note
  that the virtual fundamental class of $ D_k(\PP^N,A,B,M) $ was defined to
  be $ \frac {m^{(1)} \cdots m^{(r)}}{r!} $ times the usual one (where $ r=|M|
  $), but that on the other hand every irreducible component of the zero locus
  of $ \sigma_k $ (which is of the form $ D_k(\PP^N,A,B,M) $ for some $A$, $B$,
  $M$) gets counted $ r! $ times in the above sum, corresponding to the choice
  of order of the external components $ C^{(1)},\dots,C^{(r)} $. Hence, to
  prove the main theorem for hyperplanes in $ \PP^N $ in the case $ \sum
  \alpha = d $, we have to show that $ \sigma_k $ vanishes along $ D_k(\PP^N,
  A,B,M) $ with multiplicity $ m^{(1)} \cdots m^{(r)} $.
\end {remark}

We will now prove the main theorem for $ X=\PP^1 $ and $ Y=H $ a point, in
the case where $ \sum \alpha = d $. The proof is very similar to the proof of
\cite {V} proposition 4.8, in fact (modulo notations) identical up to the end
where the section $ \sigma_k $ comes into play, so we will only sketch these
identical parts and refer to \cite {V} for details.

\begin {proposition}[Main Theorem for $ H \subset \PP^1, \sum \alpha = d $]
  \label {main-p1-sum}
  If $ \sum \alpha = d $, then
    \[ (\alpha_k \, \psi_k + ev_k^* H) \cdot [\bar M_\alpha (\PP^1,d)]
         = [D_{\alpha,k} (\PP^1,d)]^{virt} \]
  in the Chow group of $ \bar M_\alpha (\PP^1,d) $, for all $ 1 \le k \le n $.
\end {proposition}

\begin {proof}
  Let $ D_k (\PP^1,A,B,M) $ be a component of $ D_{\alpha,k} (\PP^1,d) $. By
  equation (\ref {eqn-dim}) of definition \ref {def-d} we know that $ \sum
  \alpha^{(0)} = \sum_i m^{(i)} $, call this number $d'$. Moreover, we must
  obviously have $ r > 0 $.

  We start by defining two easier moduli spaces that model locally the
  situation at hand (in a sense that is made precise later). Fix a point $ P
  \in \PP^1 $ distinct from $H$. Let $ M \subset \bar M_{ |\alpha^{(0)}|+r}
  (\PP^1,d') $ be the closure of all degree-$ d' $ irreducible stable maps $
  (\PP^1,(x_i)_{1 \le i \le |\alpha^{(0)}|}, (y_i)_{1 \le i \le r},f) $ such
  that
    \[ f^* H = \sum_i \alpha^{(0)}_i x_i \quad \mbox {and} \quad
       f^* P = \sum_i m^{(i)} y_i. \]
  Let $ D \subset \bar M_{ |\alpha^{(0)}|+r} (\PP^1,d') $ be the closure of
  all degree-$ d' $ reducible stable maps $ (C^{(0)} \cup \cdots \cup C^{(r)},
  (x_i)_{1 \le i \le |\alpha^{(0)}|}, (y_i)_{1 \le i \le r},f) $ with $ r+1 $
  components such that
  \begin {itemize}
  \item $f$ contracts $ C^{(0)} $ to $H$, and $ C^{(i)} \cap C^{(0)} \neq
    \emptyset $ for all $ 1 \le i \le r $,
  \item $ x_i \in C^{(0)} $ for all $ 1 \le i \le |\alpha^{(0)}| $,
  \item $ (f|_{C^{(i)}})^* H = m^{(i)} (C^{(i)} \cap C^{(0)}) $ and
    $ (f|_{C^{(i)}})^* P = m^{(i)} y_i $ for all $ 1 \le i \le r $.
  \end {itemize}
  General elements of these moduli spaces look as follows (the picture
  represents the case $ \alpha = (0,4,1) $ and $ M = (2,3) $):
  \begin {center}
    \epsfig {file=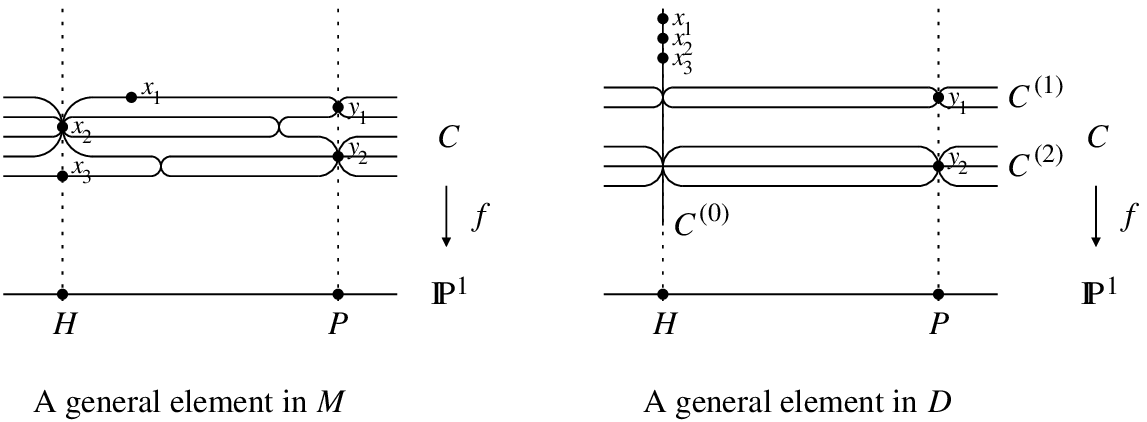}
  \end {center}
  In short, in addition to our usual multiplicity requirements for $ f^* H $
  we require multiplicities $ m^{(i)} $ over the point $P$ (so that the
  curves $ C^{(i)} $ in $D$ are ramified completely over $H$ and $P$ for $ i>0
  $).

  We are now ready to compute the multiplicity of $ \sigma_k $ to $ D_k
  (\PP^1,A,B,M) $ at a general element $ \cc' = (C',x'_1,\dots,x'_n,f') $. Let
  $ \cc = (C,(x_i),(y_i),f) $ be the unique stable map in $D$ whose internal
  component $ C^{(0)} $ is equal to the internal component of $ \cc' $, viewed
  as a marked curve whose marked points are the $ x_i $ and the points $
  C^{(0)} \cap C^{(i)} $.

  By construction, the stable maps $ \cc $ and $ \cc' $ are étale locally
  isomorphic around $ C^{(0)} $, so let $ (U,(x_i),f|_U) $ be a sufficiently
  small common étale neighborhood of $ C^{(0)} $. By \cite {V} proposition 4.3
  the deformation spaces of $ \cc $ in $M$ and $ \cc' $ in $ \bar M_\alpha
  (\PP^1,d) $ are products one of whose factors is the deformation space of $
  (U,(x_i),f|_U) $, viewed as a map from $U$ to $ \PP^1 $ satisfying the given
  multiplicity conditions at the points $ x_i $. As the section $ \sigma_k $ is
  defined on this common factor, the order of vanishing of $ \sigma_k $ along $
  D_k (\PP^1,A,B,M) $ in $ \bar M_\alpha (\PP^1,d) $ at the point $ \cc' $ is
  equal to its order of vanishing along $D$ in $M$ at the point $ \cc $.

  To simplify the calculations even further, we will now fix the marked curve $
  (C,(x_i),(y_i)) $. Consider the morphism $ \pi: M \to \bar M_{|\alpha(0)|+r}
  $ given by forgetting the map $f$ and stabilizing if necessary. Note that $
  \pi $ will contract all external components of $ \cc $ as they only have two
  special points, so $ \pi $ maps $ \cc $ to a general point of $ \bar
  M_{|\alpha(0)|+r} $. Denote by $ M' \subset M $ and $ D' \subset D $ the
  fibers of this morphism over $ \pi (\cc) $. Then the multiplicity we seek is
  equal to the multiplicity of $ \sigma_k $ along $ D' $ in $ M' $ in the point
  $ \cc $.

  But general elements in $ M' $ are actually easy to describe explicitly:
  choose $ g_1,g_2 \in \OO_{\PP^1}(d') $ with associated divisors
    \[ (g_1) = \sum_i \alpha^{(0)}_i x_i \quad \mbox {and} \quad
       (g_2) = \sum_i m^{(i)} y_i \]
  where $ x_i $ and $ y_i $ are now fixed points in $ \PP^1 $, determined
  by the element $ \pi (\cc) \in \bar M_{|\alpha(0)|+r} $. Then a general
  stable map in $ M' $ is of the form
    \[ \cc_\lambda = (\PP^1,(x_i),(y_i),f) \quad \mbox {where} \quad
       f: \PP^1 \to \PP^1, x \mapsto (\lambda g_1(x):g_2(x)) \]
  for $ \lambda \in \CC^* $. (Here we have chosen coordinates on the target
  $ \PP^1 $ such that $ H=(0:1) $ and $ P=(1:0) $.) The locus $ D' \subset M'
  $, which is set-theoretically the zero locus of $ \sigma_k $, corresponds to
  the degeneration $ \lambda \to 0 $.

  After a finite base change we can extend the family $ \{\cc_\lambda\} $
  to $ \lambda = 0 $. The central fiber $ \cc_0 $ of this extended family is
  equal to $ \cc $.

  Let $z$ be a local coordinate around $ x_k \in \PP^1 $. This means that
  $z$ is a local coordinate around $ x_k $ on all $ \cc_\lambda $ with
  $ \lambda \neq 0 $, and in fact it extends to a local coordinate around
  $ x_k $ for $ \lambda = 0 $. Consider the local trivialization of the line
  bundle $ L_k^{\otimes \alpha_k} \otimes ev_k^* \OO(H) $ given by $
  dz(x_k)^{\otimes \alpha_k} \otimes h(x_k) \mapsto 1 $ (where $ h \in H^0(
  \PP^1,\OO(H)) $ is the section vanishing at $H$ that is used to define $
  \sigma_k $). Then by construction, the section $ \sigma_k $ on the family $
  \cc_\lambda $ is given by $ \lambda \mapsto \frac {\partial^{\alpha_k}}{
  \partial z^{\alpha_k}} \lambda g_1(z)|_{z=x_k} $ in this local
  trivialization. In particular, this has a zero of first order in $ \lambda $
  at $ \lambda = 0 $. This means that the class of the zero locus of $
  \sigma_k $ on $ M' $ is
    \[ (\alpha_k \, \psi_k + ev_k^* H) \cdot [M'] = 1 \cdot [\cc_\lambda] \]
  for general $ \lambda $.

  Finally, as the automorphism group of a general $ \cc_\lambda $ is trivial,
  whereas the automorphism group of $ \cc $ is $ \ZZ_{m^{(1)}} \times \cdots
  \times \ZZ_{m^{(r)}} $, we conclude that
    \[ (\alpha_k \, \psi_k + ev_k^* H) \cdot [M'] =
       m^{(1)} \cdots m^{(r)} \cdot [\cc]. \]
  Hence the statement of the proposition follows from remark \ref {factorial}.
\end {proof}

\begin {corollary}[Main Theorem for $ H \subset \PP^N, \sum \alpha = d $]
  \label {main-pn-sum}
  If $ \sum \alpha = d $, then
    \[ (\alpha_k \, \psi_k + ev_k^* H) \cdot [\bar M_\alpha (\PP^N,d)]
         = [D_{\alpha,k} (\PP^N,d)]^{virt} \]
  in the Chow group of $ \bar M_\alpha (\PP^N,d) $, for all $ 1 \le k \le n $.
\end {corollary}

\begin {proof}
  (Compare to \cite {V} theorem 6.1.) By the previous proposition we can assume
  that $ N \ge 2 $. Consider a general element $ \cc = (C,x_1,\dots, x_n,f) $
  of a component $ D_k(\PP^N,A,B,M) $ of $ D_{\alpha,k}(\PP^N,d) $. Let $ A
  \subset H $ be a general $ (N-2) $-plane. The projection from $A$ in $ \PP^N
  $ induces a rational map $ \rho_A: \bar M_n (\PP^N,d) \dashrightarrow \bar
  M_n (\PP^1,d) $. By \cite {V} proposition 5.5 the map $ \rho_A $ is defined
  and smooth at $ \cc $. Moreover, $ \rho_A $ maps $ D_k(\PP^N,A,B,M) $ to $
  D_k(\PP^1,A,B,M) $ at the points of $ D_k(\PP^N,A,B,M) $ where it is defined,
  and the section $ \sigma_k $ on $ \bar M_\alpha (\PP^1,d) $ pulls back along
  $ \rho_A $ to the section $ \sigma_k $ on $ \bar M_\alpha (\PP^N,d) $. Hence
  the multiplicity of $ \sigma_k $ on $ \bar M_\alpha (\PP^N,d) $ along $ D_k
  (\PP^N,A,B,M) $ is the same as the multiplicity of $ \sigma_k $ on $ \bar
  M_\alpha (\PP^1,d) $ along $ D_k (\PP^1,A,B,M) $. The corollary then follows
  from proposition \ref {main-p1-sum} and remark \ref {factorial}.
\end {proof}

\begin {corollary}[Main Theorem for $ H \subset \PP^N $]
  \label {main-pn}
  We have
    \[ (\alpha_k \, \psi_k + ev_k^* H) \cdot [\bar M_\alpha (\PP^N,d)]
         = [\bar M_{\alpha+e_k} (\PP^N,d)]
           +[D_{\alpha,k} (\PP^N,d)]^{virt} \]
  in the Chow group of $ \bar M_\alpha (\PP^N,d) $, for all $ 1 \le k \le n $.
\end {corollary}

\begin {proof}
  Let $ s = d - \sum \alpha $, and let $ \alpha' = \alpha \cup (1,\dots,1) $
  such that $ \sum \alpha' = d $. By corollary \ref {main-pn-sum} we know that
  \begin {equation} \label {upstairs}
    (\alpha'_k \, \psi'_k + {ev'_k}^* H) \cdot [\bar M_{\alpha'} (\PP^N,d)]
      = [D_{\alpha',k} (\PP^N,d)]^{virt}
  \end {equation}
  for $ 1 \le k \le n $, where $ \psi'_k $ is the $k$-th cotangent line class
  on $ \bar M_{n+s} (\PP^N,d) $, and $ ev'_k $ the evaluation map $ \bar
  M_{n+s} (\PP^N,d) \to \PP^N $ at the $k$-th marked point. We will show that
  the push-forward of this equation along the morphism $ \phi: \bar M_{\alpha'}
  (\PP^N,d) \to \bar M_\alpha (\PP^N,d) $ that forgets the additional $s$
  marked points is exactly the statement of the corollary.

  First note that $ \alpha'_k = \alpha_k $ and $ ev'_k = ev_k \circ \phi $. For
  the computation of the push-forward of $ \psi'_k $ we may assume that $
  \alpha_k > 0 $, as otherwise there is no $ \psi'_k $-term in (\ref
  {upstairs}). It is well-known that $ \psi'_k = \phi^* \psi_k + \gamma $,
  where the correction term $ \gamma $ is the class of the locus of those
  stable maps $ \cc = (C,x_1,\dots,x_{n+s},f) $ where $ \phi $ contracts the
  irreducible component $Z$ of $C$ on which $ x_k $ lies, i.e.\ where $Z$ is an
  unstable component of the prestable map $ (C,x_1,\dots,x_n,f) $. This can
  only happen if $Z$ is contracted by $f$, in particular $ \sigma_k (\cc) = 0
  $, so by lemma \ref {sigma-set} the cycle $ \gamma $ must be a union of some
  of the components of $ D_k (\PP^N,A,B,M) $ of $ D_{\alpha',k} (\PP^N,d) $. To
  determine which of them occur in $ \gamma $, we can assume that $ \cc $ is a
  generic element of some $ D_k (\PP^N,A,B,M) $. It is easy to see that $ \phi
  $ contracts $ Z=C^{(0)} $ if and only if $ r=|M|=1 $, $ d^{(0)} = 0 $, and
  the marked points on $Z$ are $ x_k $ and at least one of the points $
  x_{n+1},\dots, x_{n+s} $. If there is more than one of these points on $Z$,
  the map $ \phi $ has positive-dimensional fibers on $ D_k (\PP^N,A,B,M) $,
  and hence $ \phi_* [D_k (\PP^N,A,B,M)] $ vanishes, hence we can assume that
  the marked points on $Z$ are exactly $ x_k $ and one of the forgotten points.
  Then $ \phi(\cc) $ contracts $Z$, so by remark \ref {spell-out-sum} the
  stable map $ \phi(\cc) $ will be irreducible with multiplicity $ \alpha_k+1 $
  at $ x_k $ to $H$. This means that $ \phi (D_k (\PP^N,A,B,M)) = \bar
  M_{\alpha+e_k} (\PP^N,d) $. As there is an $ s! $-fold choice of order of the
  forgotten marked points, we have shown that
    \[ \phi_* \gamma \cdot [M_{\alpha'}(\PP^N,d)]
         = s! \cdot [\bar M_{\alpha+e_k} (\PP^N,d)] \]
  and that therefore the left hand side of the push-forward of (\ref
  {upstairs}) by $ \phi $ is equal to
  \begin {equation} \label {push-lhs}
    s! \cdot (\alpha_k \, \psi_k + ev_k^* H) \cdot [\bar M_\alpha (\PP^N,d)]
      + \alpha_k s! \cdot [\bar M_{\alpha+e_k} (\PP^N,d)].
  \end {equation}
  Now we look at the right hand side of the push-forward of (\ref {upstairs})
  by $ \phi $. Consider a component $ D_k(\PP^N,A,B,M) $ of $ D_{\alpha',k}
  (\PP^N,d) $ and let $ \cc = (C,x_1,\dots,x_{n+s},f) $ be a generic element of
  this component. For the push-forward of this component by $ \phi $ to be
  non-zero, the fibers of $ \phi $ have to be zero-dimensional, i.e.\ there
  must not be a deformation of $ \cc $ inside $ D_k(\PP^N,A,B,M) $ that changes
  nothing but the position of the points $ x_{n+1},\dots,x_{n+s} $. In
  particular this means that we must have one of the following two cases:
  \begin {itemize}
  \item $ C^{(0)} $ contains none of the points $ x_{n+1},\dots,x_{n+s} $,
    i.e.\ the points $ x_{n+1},\dots,x_{n+s} $ are just the $s$ unmarked
    transverse points of intersection of $ \phi(\cc) $ with $H$. In this
    case, the map $ \phi $ does not contract any components of $C$, and it
    changes no multiplicities to $H$. Hence, the push-forward by $ \phi $ of
    all these components together is just $ s! \cdot [D_{\alpha,k} (\PP^N,d)
    ]^{virt} $.
  \item $ C^{(0)} $ is a contracted component, i.e.\ $ d^{(0)} = 0 $, $ r=|M|
    =1 $, and the marked points on $ C^{(0)} $ are exactly $ x_k $ and one of
    the points $ x_{n+1},\dots,x_{n+s} $. As above, the push-forward of such
    a component yields $ \bar M_{\alpha+e_k} (\PP^N,d) $, and it occurs with
    multiplicity $ (\alpha_k+1) \; s! $, where the factor $ \alpha_k+1 $ comes
    from the definition of the virtual fundamental class of $ D_k
    (\PP^N,A,B,M) $.
  \end {itemize}
  Put together, we have shown that the push-forward of the right hand side
  of (\ref {upstairs}) by $ \phi $ is equal to
    \[ s! \cdot [D_{\alpha,k} (\PP^N,d)]^{virt}
         + (\alpha_k+1) \; s! \cdot [\bar M_{\alpha+e_k} (\PP^N,d)]. \]
  Combining this with (\ref {push-lhs}), we get the desired result.
\end {proof}


\section {Proof of the main theorem for very ample hypersurfaces}
  \label {sec-general}

Let $X$ be a smooth complex projective variety and $Y$ a smooth very ample
hypersurface. We fix the following notation. Let $ i: Y \to X $ be the
inclusion map. For $ \beta \in H_2^+(X) $ we denote by $ \bar M_n (Y,\beta) $
the disjoint union of all moduli spaces $ \bar M_n (Y,\beta') $ for $ \beta'
\in H_2^+(Y) $ such that $ i_* \beta' = \beta $. Consider the embedding $
\varphi: X \to \PP^N $ given by the complete linear system $ |Y| $ and let $ H
\subset \PP^N $ be the hyperplane such that $ \varphi^{-1}(H) = Y $. There is
an induced morphism $ \phi: \bar M_n (X,\beta) \to \bar M_n (\PP^N,d) $, where
$ d = Y \cdot \beta $. In this section we will show that the ``pull-back'' of
the main theorem for $ H \subset \PP^N $ by $ \phi $ yields the main theorem
for $ Y \subset X $. The most difficult part of the proof is to show that
the spaces $ D_{\alpha,k} (\PP^N,d) $ pull back to $ D_{\alpha,k} (X,\beta) $
(proposition \ref {pullback-d}). Recall that curves in $ D_{\alpha,k} (X,\beta)
$ are reducible curves with one component in $Y$ (and some multiplicity
conditions). Hence we will show first that the moduli spaces of curves in $Y$
(lemma \ref {virt-comp}) and those of reducible curves in $X$ (lemma \ref
{comb-comp}) pull back nicely under $ \phi $.

\begin {convention}
  In this section, all occurring spaces are equipped with virtual fundamental
  classes as follows.
  \begin {itemize} \itemsep 0mm
  \item The moduli spaces of stable maps $ \bar M_n (\cdot,\cdot) $ have
    virtual fundamental classes constructed e.g.\ in \cite {B}, \cite {BF}.
  \item The moduli spaces $ \bar M_\alpha (\cdot,\cdot) $, $ D_k (\dots) $,
    and $ D_{\alpha,k}(\dots) $ have virtual fundamental classes constructed
    in definitions \ref {def-fund}, \ref {def-dc}, and \ref {def-d},
    respectively.
  \item The varieties $Y$, $X$, $H$, and $ \PP^N $ are equipped with their
    usual fundamental class.
  \item The virtual fundamental class of a disjoint union of spaces is the
    sum of the virtual fundamental classes of its components.
  \item In any fiber product $ V_1 \times_V V_2 $ occurring in this section,
    $V$ will always be smooth and equipped with the usual fundamental class.
    The virtual fundamental class of the fiber product is then taken to be
    the one induced by the virtual fundamental classes of $ V_1 $ and $ V_2 $
    in the sense of remark \ref {ind-fund}.
  \end {itemize}
  When we say that two spaces $ V_1 $ and $ V_2 $ are equal we will always
  mean that $ V_1 $ and $ V_2 $ are isomorphic and that $ [V_1]^{virt} =
  [V_2]^{virt} $ under this isomorphism. We will write this as $ V_1 \equiv
  V_2 $.
\end {convention}

\begin {lemma} \label {virt-comp}
  For any $ n \ge 0 $ and $ \beta \in H_2^+(X) $ we have
    \[ \bar M_n (Y,\beta) \equiv
       \bar M_n (H,d) \times_{\bar M_n (\PP^N,d)} \bar M_n (X,\beta). \]
\end {lemma}

\begin {proof}
  As $ Y = H \cap X \subset \PP^N $, it follows from the definitions that
  the diagram of inclusions
  \begin {equation} \label {cart-mod} \begin {split} \xymatrix @M=5pt {
    \bar M_n (Y,\beta) \ar[r] \ar[d] & \bar M_n (X,\beta) \ar[d] \\
    \bar M_n (H,d) \ar[r]^-\psi        & \bar M_n (\PP^N,d)
  } \end {split} \end {equation}
  is cartesian. We denote by $ \pi_X: \bar M_{n+1} (X,\beta) \to \bar M_n
  (X,\beta) $ the universal curve and by $ f_X: \bar M_{n+1} (X,\beta) \to X $
  its evaluation map, and similarly for the moduli spaces of maps to $Y$, $H$,
  and $ \PP^N $. Applying the functor $ R{\pi_Y}_* f_Y^* $ to the distinguished
  triangle
  \begin {equation} \label {cot-y}
    L_X|_Y \to L_Y \to L_{Y/X} \to L_X|_Y[1]
  \end {equation}
  on $Y$, we get the distinguished triangle
  \begin {align*}
    R&{\pi_Y}_* (f_X^* L_X)|_{\bar M_{n+1} (Y,\beta)} \to
    R{\pi_Y}_* f_Y^* L_Y \to
    R{\pi_Y}_* (f_H^* L_{H/\PP^N})|_{\bar M_{n+1} (Y,\beta)} \\
    &\to R{\pi_Y}_* (f_X^* L_X)|_{\bar M_{n+1} (Y,\beta)} [1]
  \end {align*}
  on $ \bar M_n (Y,\beta) $. By \cite {B} proposition 5, the vector bundle $
  f_X^* L_X $ is quasi-isomorphic to a complex $K$ of vector bundles on $ \bar
  M_{n+1} (X,\beta) $ such that $ R{\pi_X}_* K $ is also a complex of vector
  bundles. As $ \pi_X $ is flat, it follows from the theorem on cohomology and
  base change that $ (R{\pi_X}_* K)_{\bar M_n (Y,\beta)} = R{\pi_Y}_* (K|_{\bar
  M_{n+1} (Y,\beta)}) $. The same argument applies to $ f_H^* L_{H/\PP^N} $
  instead of $ f_X^* L_X $, so we arrive at the distinguished triangle
  \begin {equation} \label {dist-y} \begin {split}
    (R&{\pi_X}_* f_X^* L_X)|_{\bar M_n (Y,\beta)} \to
    R{\pi_Y}_* f_Y^* L_Y \to
    (R{\pi_H}_* f_H^* L_{H/\PP^N})|_{\bar M_n (Y,\beta)} \\
    &\to (R{\pi_X}_* f_X^* L_X)|_{\bar M_n (Y,\beta)} [1].
  \end {split} \end {equation}
  Starting with the distinguished triangle of $ L_{H/\PP^N} $ instead of
  $ L_{Y/X} $ in (\ref {cot-y}), the same calculation as above shows that
  we also have a distinguished triangle on $ \bar M_n (H,d) $
  \begin {equation*} \begin {split}
    (R&{\pi_{\PP^N}}_* f_{\PP^N}^* L_{\PP^N})|_{\bar M_n (H,d)} \to
    R{\pi_H}_* f_H^* L_H \to
    R{\pi_H}_* f_H^* L_{H/\PP^N} \\
    &\to (R{\pi_{\PP^N}}_* f_{\PP^N}^* L_{\PP^N})|_{\bar M_n (H,d)} [1].
  \end {split} \end {equation*}
  But the first and second term in this sequence are just $ L_{\bar M_n
  (\PP^N,d)/\MM_n}|_{\bar M_n (H,d)} $ and $ L_{\bar M_n (H,d)/\MM_n} $, where
  $ \MM_n $ denotes the stack of prestable $n$-pointed rational curves. Hence
  we see that $ R{\pi_H}_* f_H^* L_{H/\PP^N} = L_{\bar M_n (H,d) / \bar M_n
  (\PP^N,d)} $. So (\ref {dist-y}) becomes
  \begin {equation*} \begin {split}
    (R&{\pi_X}_* f_X^* L_X)|_{\bar M_n (Y,\beta)} \to
    R{\pi_Y}_* f_Y^* L_Y \to
    {L_{\bar M_n (H,d) / \bar M_n (\PP^N,d)}}|_{\bar M_n (Y,\beta)} \\
    &\to (R{\pi_X}_* f_X^* L_X)|_{\bar M_n (Y,\beta)} [1].
  \end {split} \end {equation*}
  As the first two terms in this sequence are the relative obstruction theories
  of $ \bar M_n (X,\beta) $ and $ \bar M_n (Y,\beta) $ over $ \MM_n $,
  respectively, we get a homomorphism of this distinguished triangle to
  \begin {equation*} \begin {split}
    &{L_{\bar M_n (X,\beta) / \MM_n}}|_{\bar M_n (Y,\beta)} \to
    L_{\bar M_n (Y,\beta) / \MM_n} \to
    L_{\bar M_n (Y,\beta) / \bar M_n (X,\beta)} \\
    &\;\;\to {L_{\bar M_n (X,\beta) / \MM_n}}|_{\bar M_n (Y,\beta)} [1].
  \end {split} \end {equation*}
  Hence, by \cite {BF} proposition 7.5 it follows that $ \psi^! [\bar M_n
  (X,\beta)]^{virt} = [\bar M_n (Y,\beta)]^{virt} $ in (\ref {cart-mod}).
  This proves the lemma.
\end {proof}

\begin {lemma} \label {comb-comp}
  Let $ n^{(i)} \ge 0 $ and $ d^{(i)} \ge 0 $ such that $ \sum_i n^{(i)} = n $
  and $ \sum_i d^{(i)} = d $. Then
  \begin {align*}
    \coprod_{(\beta^{(i)})} &\left(
      \bar M_{n^{(0)}+r} (X,\beta^{(0)}) \times_{X^r}
      \prod_{i=1}^r \bar M_{n^{(i)}+1} (X,\beta^{(i)})
    \right)
    \equiv \\
    &\kern-3mm \left(
      \bar M_{n^{(0)}+r} (\PP^N,d^{(0)}) \times_{(\PP^N)^r}
      \prod_{i=1}^r \bar M_{n^{(i)}+1} (\PP^N,d^{(i)})
    \right) \times_{\bar M_n (\PP^N,d)} \bar M_n (X,\beta),
  \end {align*}
  where the union is taken over all $ (\beta^{(i)}) $ with $ Y \cdot \beta^{
  (i)} = d^{(i)} $ for all $i$, and where the maps to $ X^r $ and $ (\PP^N)^r
  $ are given in the same way as in definition \ref {def-dc}.
\end {lemma}

\begin {proof}
  In the language of \cite {BM}, let $ \tau $ be the graph corresponding to
  rational curves with components $ C^{(0)},\dots,C^{(r)} $ such that $ C^{(0)}
  \cap C^{(i)} \neq\emptyset $ for all $ i > 0 $ and $ C^{(i)} $ has $ n^{(i)}
  $ marked points for $ i \ge 0 $. Let $ \MM_n $ be the stack of prestable
  $n$-pointed rational curves, and let $ \MM_\tau \subset \MM_n $ be the
  substack of $ \tau $-marked prestable curves, as defined in \cite {BM}
  definition 2.6. Moreover, we will abbreviate the moduli spaces in the large
  brackets in the statement of the lemma as $ \bar M_\tau (X,(\beta^{(i)}))
  $ and $ \bar M_\tau (\PP^N,(d^{(i)})) $, respectively.

  Consider the commutative diagram
  \xydiag {
    \bar M_\tau (X,(\beta^{(i)})) \ar[r] \ar[d] &
      \bar M_\tau (\PP^N,(d^{(i)})) \ar[r] \ar[d] &
      \MM_\tau \ar[d]^{\psi} \\
    \bar M_n (X,\beta) \ar[r] &
      \bar M_n (\PP^N,d) \ar[r] &
      \MM_n
  }
  where none of the maps involves stabilization of the underlying prestable
  curves. By \cite {B} lemma 10, the right square and the big square are
  cartesian, so the left one is also cartesian. Moreover, by the same lemma, $
  \psi^! [\bar M_n (X,\beta)]^{virt} = [\bar M_\tau (X,(\beta^{(i)}))]^{virt}
  $.
\end {proof}

\begin {proposition} \label {pullback-d}
  For any $ 1 \le k \le n $ we have
    \[ D_{\alpha,k}(X,\beta) \equiv
         D_{\alpha,k}(\PP^N,d)
         \times_{\bar M_n (\PP^N,d)} \bar M_n (X,\beta). \]
  In particular, the moduli spaces $ D_k (X,A,B,M) $ satisfying equation
  (\ref {eqn-dim}) of definition \ref {def-d} are proper substacks of $ \bar
  M_\alpha (X,\beta) $ of expected codimension one.
\end {proposition}

\begin {proof}
  We consider a component $ D_k (\PP^N,A,(d^{(i)}),M) $ of $ D_{\alpha,k}
  (\PP^N,d) $ and show that the fiber product of this component with $ \bar M_n
  (X,\beta) $ over $ \bar M_n (\PP^N,d) $ is the union of all $ D_k
  (X,A,(\beta^{(i)}),M) $ such that $ Y \cdot \beta^{(i)} = d^{(i)} $.

  We start with the pull-back compatibility statement for general curves of the
  form $ C^{(0)} \cup \cdots \cup C^{(r)} $ with $ C^{(0)} \cap C^{(i)} \neq
  \emptyset $, as given in lemma \ref {comb-comp}. Taking the fiber product
  of this equation with $ \bar M_{n^{(0)}+r} (H,d^{(0)}) $ over $ \bar M_{
  n^{(0)}+r} (\PP^N,d^{(0)}) $ (i.e.\ requiring the central component $
  C^{(0)} $ to lie in $H$) and using lemma \ref {virt-comp} on the left hand
  side yields
  \begin {align*}
    \coprod_{(\beta^{(i)})} &\left(
      \bar M_{n^{(0)}+r} (Y,\beta^{(0)})
      \times_{X^r} \prod_{i=1}^r \bar M_{n^{(i)}+1} (X,\beta^{(i)})
    \right)
    \equiv \\
    &\kern-3mm \left(
      \bar M_{n^{(0)}+r} (H,d^{(0)}) \times_{(\PP^N)^r}
      \prod_{i=1}^r \bar M_{n^{(i)}+1} (\PP^N,d^{(i)})
    \right) \times_{\bar M_n (\PP^N,d)} \bar M_n (X,\beta).
  \end {align*}
  This can obviously be written in a more complicated way as
  {\small \begin {align*}
    \coprod_{(\beta^{(i)})} &\left(
      \bar M_{n^{(0)}+r} (Y,\beta^{(0)})
      \times_{Y^r} \left(
        H^r \times_{(\PP^N)^r}
        \prod_{i=1}^r \bar M_{n^{(i)}+1} (X,\beta^{(i)})
      \right)
    \right)
    \equiv \\
    &\kern-3mm \left(
      \bar M_{n^{(0)}+r} (H,d^{(0)}) \times_{H^r} \left(
        H^r \times_{(\PP^N)^r}
      \prod_{i=1}^r \bar M_{n^{(i)}+1} (\PP^N,d^{(i)}) \right)
    \right) \times_{\bar M_n (\PP^N,d)} \bar M_n (X,\beta).
  \end {align*}}%
  Note that $ H \times_{\PP^N} \bar M_{n^{(i)}+1} (\PP^N,d^{(i)}) \equiv
  \bar M_{\tilde \alpha^{(i)}} (\PP^N,d^{(i)}) $ for all $ i>0 $, where $
  \tilde \alpha^{(i)} = (0,\dots,0,1) $. So we get
  {\small \begin {align*}
    \coprod_{(\beta^{(i)})} &\left(
      \bar M_{n^{(0)}+r} (Y,\beta^{(0)})
      \times_{Y^r}
        \prod_{i=1}^r
          \bar M_{\tilde \alpha^{(i)}} (\PP^N,d^{(i)})
            \times_{\bar M_{n^{(i)}+1} (\PP^N,d^{(i)})}
          \bar M_{n^{(i)}+1} (X,\beta^{(i)})
    \right)
    \equiv \\
    &\kern-3mm \left(
      \bar M_{n^{(0)}+r} (H,d^{(0)}) \times_{H^r}
      \prod_{i=1}^r \bar M_{\tilde \alpha^{(i)}} (\PP^N,d^{(i)})
    \right) \times_{\bar M_n (\PP^N,d)} \bar M_n (X,\beta).
  \end {align*}}%
  Finally, we take the fiber product of this equation with $ \bar M_{
  \alpha^{(i)} \cup (m^{(i)})} (\PP^N,d) $ over $ \bar M_{\tilde \alpha^{(i)}}
  (\PP^N,d) $ for all $ i>0 $, yielding the same equation with the $ \tilde
  \alpha^{(i)} $ replaced by $ \alpha^{(i)} \cup (m^{(i)}) $. By definition,
  this is then exactly the equation stated in the proposition.
\end {proof}

We are now ready to give the proof of our main theorem.

\begin {proof}[Proof (of theorem \ref {main-thm})]
  Consider the cartesian diagram
  \xydiag {
    \bar M_\alpha (X,\beta) \ar[r] \ar[d] &
      \bar M_\alpha (\PP^N,d) \ar[d] \\
    \bar M_n (X,\beta) \ar[r]^-\phi &
      \bar M_n (\PP^N,d).
  }
  The main theorem for $ H \subset \PP^N $ (see corollary \ref {main-pn})
  gives an equation in the Chow group of $ \bar M_\alpha (\PP^N,d) $. We pull
  this equation back by $ \phi $ to get an equation in the Chow group of
  $ \bar M_\alpha (X,\beta) $. As the morphism $ \phi $ does not involve any
  contractions of the underlying prestable curves, the cotangent line class $
  \psi_k $ on $ \bar M_n (\PP^N,d) $ pulls back to the cotangent line class $
  \psi_k $ on $ \bar M_n (X,\beta) $. So by definition the left hand side of
  corollary \ref {main-pn} pulls back to $ (\alpha_k \psi_k + ev^* Y) \cdot
  [\bar M_\alpha (X,\beta)]^{virt} $. In the same way, $ [\bar M_{\alpha+e_k}
  (\PP^N,d)] $ pulls back to $ [\bar M_{\alpha+e_k} (X,\beta)]^{virt} $.
  Finally, proposition \ref {pullback-d} shows that $ [D_{\alpha,k}
  (\PP^N,d)]^{virt} $ pulls back to $ [D_{\alpha,k} (X,\beta)]^{virt} $.
\end {proof}

\begin {remark}
  We expect that the statement of the main theorem \ref {main-thm} is true even
  under weaker assumptions on the hypersurface $Y$. For example, if $Y$ is not
  very ample but the complete linear system $ |Y| $ on $X$ is base-point free,
  we still get a morphism $ X \to \PP^N $ defined by $ |Y| $. The definition of
  the moduli spaces of relative invariants essentially carries over without
  change to this case. The main (but probably little) problem is that the
  morphism $ \phi $ in the cartesian diagram of definition \ref {def-fund} now
  may involve stabilization of the underlying prestable curves.  This makes
  many points in the arguments of this paper more subtle, but we expect that a
  version of the main theorem can be proven also in this case.
\end {remark}


\section {Enumerative applications} \label {sec-enum}

As usual, the first thing to do to get enumerative results from moduli spaces
of maps is to define invariants by intersecting the virtual fundamental class
of the moduli space with various cotangent line classes and pull-backs of
classes via evaluation maps. Note that from the spaces $ \bar M_\alpha
(X,\beta) $ we always have evaluation maps $ ev_k $ to $X$ for $ 1 \le k \le
|\alpha| $, and in addition evaluation maps $ \tilde {ev}_k $ to $Y$ for all
$k$ with $ \alpha_k > 0 $.

\begin {definition}
  Let $ \beta \in H_2^+(X) $, $ n \ge 0 $, $ k_1,\dots,k_n \ge 0 $, and $
  \gamma_1,\dots,\gamma_n \in A^*(X) $. Then the restricted Gromov-Witten
  invariants of $Y$ are defined as
    \[ I_{n,\beta}^Y (\gamma_1 \psi^{k_1},\dots,\gamma_n \psi^{k_n})
         = ev_1^* \gamma_1 \cdot \psi_1^{k_1} \cdots
           ev_n^* \gamma_n \cdot \psi_n^{k_n}
           \cdot [\bar M_n (Y,\beta)]^{virt} \in \QQ \]
  if $ \sum_i (\codim \gamma_i + k_i) = \vdim \bar M_n (Y,\beta) $. Similarly,
  for any $ \alpha = (\alpha_1,\dots,\alpha_n) $ the restricted relative
  Gromov-Witten invariants of $ Y \subset X $ are defined as
    \[ I_{\alpha,\beta} (\gamma_1 \psi^{k_1},\dots,\gamma_n \psi^{k_n})
         = ev_1^* \gamma_1 \cdot \psi_1^{k_1} \cdots
           ev_n^* \gamma_n \cdot \psi_n^{k_n}
           \cdot [\bar M_\alpha (X,\beta)]^{virt} \in \QQ \]
  if $ \sum_i (\codim \gamma_i + k_i) = \vdim \bar M_\alpha (X,\beta) $. This
  definition can obviously be generalized in the following two ways:
  \begin {enumerate}
  \item \label {restr-coh}
    We can take cohomology classes $ \tilde \gamma_k \in A^*(Y) $ and the
    evaluation maps $ \tilde {ev}_k $ to $Y$, instead of $ \gamma_k \in A^*(X)
    $ and $ ev_k $ (provided that $ \alpha_k > 0 $ in the case of the relative
    invariants). We will apply the same notation in this case and just mark the
    cohomology classes that are pulled back from $Y$ by a tilde.
  \item \label {restr-hom}
    For the absolute invariants, we could use a homology class on $Y$ instead
    of summing over all homology classes on $Y$ that push forward to a given
    class on $X$. (We will never do this in this paper, however.)
  \end {enumerate}
  The invariants obtained in this way are called the (unrestricted)
  Gromov-Witten invariants of $Y$, or relative Gromov-Witten invariants of $
  Y \subset X $, respectively.
\end {definition}

\begin {remark} \label {no-restr}
  Often the restricted invariants are really not restricted at all. As for
  generalization \ref {restr-coh} in the above definition, in many cases every
  algebraic cohomology class in $Y$ comes from a (rational) algebraic
  cohomology class in $X$, notably if the dimension of $Y$ is odd (by the
  Lefschetz theorem) or if $ X=\PP^N $ and $Y$ is a generic hypersurface that
  is not a quadric or the cubic surface (by \cite {S} proposition 2.1). Again
  by the Lefschetz theorem, \ref {restr-hom} is no generalization if the
  dimension of $Y$ is at least 3.
\end {remark}

\begin {remark} \label {diagonal-trick}
  If we intersect the main theorem \ref {main-thm}
    \[ (\alpha_k \, \psi_k + ev_k^* Y) \cdot [\bar M_\alpha (X,\beta)]^{virt}
         = [\bar M_{\alpha+e_k} (X,\beta)]^{virt}
           + [D_{\alpha,k} (X,\beta)]^{virt} \]
  with suitably many cotangent line classes or pull-backs from classes on $X$
  or $Y$ by the evaluation maps, we obviously get many relations among the
  relative Gromov-Witten invariants of $ Y \subset X $, the Gromov-Witten
  invariants of $X$ (for $ \alpha=(0,\dots,0) $), and the Gromov-Witten
  invariants of $Y$ (as the moduli spaces of stable maps to $Y$ are included as
  factors in the spaces $ D_{\alpha,k} (X,\beta) $). As for $ D_{\alpha,k}
  (X,\beta) $ one uses the usual ``diagonal trick'' to express a component
    \[ D_k(X,A,B,M) = \bar M_{|\alpha^{(0)}|+r} (Y,\beta^{(0)}) \times_{Y^r}
         \prod_{i=1}^r \bar M_{\alpha^{(i)} \cup (m^{(i)})} (X,\beta^{(i)}) \]
  (and its virtual fundamental class) by the cartesian diagram
  \xydiag {
    D_k (X,A,B,M) \ar[r] \ar[d] &
      \bar M_{|\alpha^{(0)}|+r} (Y,\beta^{(0)}) \times
      \prod_{i=1}^r \bar M_{\alpha^{(i)} \cup (m^{(i)})} (X,\beta^{(i)})
      \ar[d]^{ev} \\
    Y^r \ar[r]^-{\Delta^r} & Y^r \times Y^r,
  }
  i.e.\ intersection products on $ D_k (X,A,B,M) $ become intersection products
  of the same classes on products of moduli spaces of (absolute and relative)
  stable maps, with additional classes coming from the diagonal. So the term $
  [D_{\alpha,k} (X,\beta)]^{virt} $ in the main theorem will turn into a sum of
  products of Gromov-Witten invariants of $Y$ and relative Gromov-Witten
  invariants of $ Y \subset X $.
\end {remark}

\begin {remark}
  In what follows we only want to look at the restricted (relative)
  Gromov-Witten invariants. It is not obvious that this is possible, as even if
  we only use pull-backs of classes from $X$ at the marked points $
  x_1,\dots,x_n $, the classes from the diagonal trick in the terms $
  D_{\alpha,k} (X,\beta) $ (see above) will throw in classes from $Y$. To see
  that these do not do any harm we will first show in the next two lemmas that
  absolute as well as relative invariants vanish if they contain exactly one
  class from $Y$ and this class lies in the orthogonal complement $
  A^*(X)^\perp $ of $ i^* A^*(X) $ in $ A^*(Y) $. (These lemmas can obviously
  be skipped if $ A^*(X)^\perp = \emptyset $, which is often the case by remark
  \ref {no-restr}).
\end {remark}

\begin {lemma} \label {gw-vanish}
  Let $ \tilde \gamma_1 \in A^*(X)^\perp $ and $ \gamma_2,\dots,\gamma_n \in
  A^*(X) $. Then for any $ \beta \in H_2^+(X) $ we have $ I_{n,\beta}^Y (\tilde
  \gamma_1 \psi^{k_1},\gamma_2 \psi^{k_2},\dots,\gamma_n \psi^{k_n}) = 0 $.
\end {lemma}

\begin {proof}
  (This is a variant of proposition 4 in \cite {P}.) Consider the cartesian
  diagram (see lemma \ref {virt-comp})
  \xydiag {
    & Y \ar[r]^-i & X \\
    \bar M_n (Y,\beta) \ar[r] \ar[d] &
      \bar M_{\tilde \alpha} (X,\beta) \ar[r] \ar[d] \ar[u]^{\tilde {ev}_1} &
      \bar M_n (X,\beta) \ar[d]^\phi \ar[u]^{ev_1} \\
    \bar M_n (H,d) \ar[r] &
      \bar M_{\tilde \alpha} (H,d) \ar[r]^-j &
      \bar M_n (\PP^N,d)
  }
  where $ \tilde \alpha = (1,0,\dots,0) $. Let $ \pi: \bar M_{n+1} (\PP^N,d)
  \to \bar M_n (\PP^N,d) $ be the universal map and $ f: \bar M_{n+1} (\PP^N,
  d) \to \PP^N $ its evaluation map. Let $E$ be the kernel of the surjective
  bundle morphism $ \pi_* f^* \OO(H) \to ev_1^* \OO(H) $ given by evaluation.
  By \cite {P} construction 2.1 and proposition 4 we have that $ [\bar M_n
  (H,d)] = j^* (c_{top}(E) \cdot [\bar M_n(\PP^N,d)]) $. Intersecting with
  $ [\bar M_n (X,\beta)]^{virt} $ yields by lemma \ref {virt-comp}
    \[ [\bar M_n (Y,\beta)]^{virt} =
         i^! (\phi^* c_{top} (E) \cdot [\bar M_n (X,\beta)]^{virt}) \]
  on $ \bar M_{\tilde \alpha} (X,\beta) $. Moreover, the class $ \gamma =
  \psi^{k_1} \cdot ev_2^* \gamma_2 \cdot \psi^{k_2} \cdots ev_n^* \gamma_n
  \cdot \psi^{k_n} $ is actually defined on $ \bar M_n (X,\beta) $. Therefore
  we get
  \begin {align*}
    I_{n,\beta}^Y (\tilde \gamma_1 \psi^{k_1},
      \gamma_2 \psi^{k_2},\dots,\gamma_n \psi^{k_n})
    &= \tilde \gamma_1 \cdot \tilde {ev}_{1*} i^!
      (\gamma \cdot \phi^* c_{top}(E) \cdot [\bar M_n (X,\beta)]^{virt}) \\
    &= \tilde \gamma_1 \cdot i^* ev_{1*}
      (\gamma \cdot \phi^* c_{top}(E) \cdot [\bar M_n (X,\beta)]^{virt}) \\
    &=0
  \end {align*}
  as $ \tilde \gamma_1 \in A^*(X)^\perp $.
\end {proof}

\begin {lemma} \label {rel-vanish}
  Assume that $ \alpha_1 > 0 $. Let $ \tilde \gamma_1 \in A^*(X)^\perp $ and
  $ \gamma_2,\dots,\gamma_n \in A^*(X) $. Then $ I_{\alpha,\beta} (\tilde
  \gamma_1 \psi^{k_1},\gamma_2 \psi^{k_2},\dots,\gamma_n \psi^{k_n}) = 0 $.
\end {lemma}

\begin {proof}
  We prove the statement by induction on $ d = Y \cdot \beta $, $n$, and $ \sum
  \alpha $, in that order. This means: if we want to prove the statement for
  an invariant with certain values of $d$, $n$, and $ \sum \alpha $, we assume
  that it is true for all invariants having
  \begin {enumerate}
  \item smaller $d$, or
  \item the same $d$ and smaller $n$, or
  \item the same $d$, the same $n$, and smaller $ \sum \alpha $.
  \end {enumerate}
  For $ \sum \alpha = 1 $, i.e.\ $ \alpha = (1,0,\dots,0) $, the statement
  follows by exactly the same calculation as in the proof of lemma \ref
  {gw-vanish}, just leaving out the factor $ c_{top}(E) $. So we can assume
  that $ \sum \alpha > 1 $. If $ \alpha_1 > 1 $ set $ k=1 $, otherwise choose
  any $ k>1 $ with $ \alpha_k > 0 $. By the main theorem
  \ref {main-thm} we have
    \[ ((\alpha_k-1) \, \psi_k + ev_k^* Y)
         \cdot [\bar M_{\alpha-e_k} (X,\beta)]^{virt}
         = [\bar M_\alpha (X,\beta)]^{virt}
           + [D_{\alpha-e_k,k} (X,\beta)]^{virt}. \]
  Intersect this equation with $ \tilde {ev}_1^* \tilde \gamma_1 \cdot
  \psi^{k_1} \cdot ev_2^* \gamma_2 \cdot \psi^{k_2} \cdots ev_n^* \gamma_n
  \cdot \psi^{k_n} $. The first term on the right hand side is then exactly
  the desired invariant. We will show that all other terms vanish.

  The term on the left hand side has the same $d$ and $n$, and smaller $ \sum
  \alpha $. The invariant coming from the $ \psi_k $-summand has exactly
  one class in $ A^*(X)^\perp $ and hence vanishes by the induction
  hypothesis. The same is true for the invariant coming from the $ ev_k^* Y
  $-term if $ k>1 $. If $ k=1 $, all classes in the invariant come from $X$,
  but the invariant contains the class $ ev_1^* Y \cdot \tilde {ev}_1^*
  \tilde \gamma_1 = \tilde {ev}_1^* (\tilde \gamma_1 \cdot i^* Y) $, which is
  zero as $ \tilde \gamma_1 \in A^*(X)^\perp $. Hence the left hand side of
  the equation vanishes.

  Now we look at the terms $ D_k (X,A,B,M) $ on the right hand side that
  give products of (relative) invariants by the diagonal trick as described 
  in remark \ref {diagonal-trick}. Note that the class of the diagonal in
  $ Y \times Y $ is $ \sum_i T_i \otimes T_i^\vee $, where $ \{T_i\} $
  is a basis of $ A^*(Y) $. If we choose this basis such that it respects the
  orthogonal decomposition $ A^*(Y) = i^* A^*(X) \oplus A^*(X)^\perp $, then
  $ T_i \in A^*(X)^\perp $ if and only if $ T_i^\vee \in A^*(X)^\perp $. Hence
  the $i$-th diagonal (where $ 1 \le i \le r $) will contribute one class
  each to the invariants for $ C^{(0)} $ and $ C^{(i)} $, and either both
  of them are in $ A^*(X)^\perp $ or none of them.

  For a given term $ D_k (X,A,B,M) $, the components $ C^{(i)} $ for $ i>0 $
  all have either smaller $d$, or the same $d$ and smaller $n$ (the latter
  happens only if $ r=1 $ and $ \beta^{(0)} = 0 $). Hence by induction
  hypothesis ($ i>0 $) or lemma \ref {gw-vanish} ($ i=0 $), we know for any $ i
  \ge 0 $ that the invariant for $ C^{(i)} $ vanishes if it contains exactly
  one class from $ A^*(X)^\perp $. We show that this has always to be the case
  for at least one $i$. Assume that this is not true. We distinguish two cases:
  \begin {enumerate}
  \item $ x_1 \in C^{(0)} $. Then the external components $ C^{(i)} $ can
    have at most one class from $ A^*(X)^\perp $, namely the class from the
    diagonal. Hence by our assumption, they have no such class, i.e.\ the
    diagonal contributes a class from $ i^* A^*(X) $ to $ C^{(i)} $ and hence
    also to $ C^{(0)} $. But then the invariant for $ C^{(0)} $ has exactly
    one class from $ A^*(X)^\perp $, namely $ \tilde \gamma_1 $, which is
    a contradiction.
  \item $ x_1 \in C^{(i)} $ for some $ i>0 $. Then by our assumption, the
    diagonals must contribute a class from $ A^*(X)^\perp $ to $ C^{(i)} $,
    and a class from $ i^* A^*(X) $ to all other $ C^{(j)} $ with $ j>0 $.
    But then we have again exactly one class from $ A^*(X)^\perp $ in $ C^{(0)}
    $, namely the one from the $i$-th diagonal. This is again a contradiction.
  \end {enumerate}
  This shows the lemma.
\end {proof}

\begin {corollary} \label {gw-hyp}
  Let $X$ be a smooth projective variety and $ Y \subset X $ a smooth very
  ample hypersurface. Assume that the Gromov-Witten invariants of $X$ are
  known. Then there is an explicit algorithm to compute the restricted
  Gromov-Witten invariants of $Y$ as well as the restricted relative
  Gromov-Witten invariants of $ Y \subset X $.
\end {corollary}

\begin {proof}
  This is now straightforward. We will compute the absolute and relative
  invariants at the same time, and we will use recursion on the same variables
  as in the previous lemma.

  Assume that we want to compute a relative invariant $ I_{\alpha, \beta}
  (\gamma_1 \psi^{k_1},\dots,\gamma_n \psi^{k_n}) $.  If $ \sum \alpha = 0 $
  then this is a Gromov-Witten invariant on $X$ and therefore assumed to be
  known. So we can assume that $ \sum \alpha > 0 $. On the other hand, we can
  also assume that $ \sum \alpha \le Y \cdot \beta = d $, as otherwise the
  invariant is zero anyway by definition.

  Choose $k$ such that $ \alpha_k > 0 $ and intersect the main theorem \ref
  {main-thm}
  \begin {equation} \label {main-thm-shift}
    ((\alpha_k-1) \, \psi_k + ev_k^* Y)
      \cdot [\bar M_{\alpha-e_k} (X,\beta)]^{virt}
    = [\bar M_\alpha (X,\beta)]^{virt}
      + [D_{\alpha-e_k,k} (X,\beta)]^{virt}
  \end {equation}
  with $ ev_1^* \gamma_1 \cdot \psi^{k_1} \cdots ev_n^* \gamma_n \cdot \psi^{
  k_n} $. Then the first term on the right hand side is the invariant that we
  want to compute. We will show that all other terms in the equation are
  recursively known.

  This is obvious for the invariants on the left hand side, since they have the
  same $d$, same $n$, and smaller $ \sum \alpha $. Now look at a term coming
  from $ D_k (X,A,B,M) $ on the right hand side, it is a product of invariants
  for the components $ C^{(i)} $ for $ i=0,\dots,r $. First we will show that
  we only get products of \emph {restricted} invariants. The invariant for the
  components $ C^{(i)} $ for $ i>0 $ can have at most one class from $
  A^*(X)^\perp $, namely from the diagonal. But if it has exactly one it
  vanishes by lemma \ref {rel-vanish}, so it has none. This means that it is a
  restricted invariant, and moreover that the diagonal contributes only classes
  from $ A^*(X) $ to the invariant for $ C^{(0)} $. This means that the
  invariant for $ C^{(0)} $ is also a restricted one.

  Now, as in the previous lemma, the invariants for the components $ C^{(i)} $
  for $ i>0 $ all have either smaller $d$, or the same $d$ and smaller $n$, and
  are therefore recursively known. The Gromov-Witten invariant for the
  component $ C^{(0)} $ can certainly have no bigger $d$. We will show now
  that it cannot have the same $d$ either. Assume the contrary, then we must
  have $ r=0 $. But then the dimension condition says
  \begin {align} \label {dim-check} \begin {split}
    \vdim \bar M_\alpha (X,\beta) &= \vdim \bar M_n (Y,\beta) \\
    \iff \vdim \bar M_n (X,\beta) - \sum \alpha &= \vdim \bar M_n (X,\beta)
      -d-1,
  \end {split} \end {align}
  i.e.\ $ \sum \alpha = d+1 > d $, which is a contradiction. Hence also the
  invariant for $ C^{(0)} $ has smaller $d$. In summary, we have seen that we
  can compute the desired relative Gromov-Witten invariant.

  Now we compute the absolute Gromov-Witten invariants for the same values of
  $d$ and $n$. Assume that there is such an invariant $ I_{n,\beta}^Y (\gamma_1
  \psi^{k_1},\gamma_2 \psi^{k_2},\dots,\gamma_n \psi^{k_n}) $. Without loss
  of generality we may assume that $ n>0 $ (if $ n=0 $ we can just add one
  marked point and require it to be on $Y$, which changes the invariant only
  by a factor of $d$ according to the divisor axiom). Set $ \alpha = (d+1,0,
  \dots,0) $. Now consider exactly the same equation (\ref {main-thm-shift})
  as above and intersect it again with $ ev_1^* \gamma_1 \cdot \psi^{k_1}
  \cdots ev_n^* \gamma_n \cdot \psi^{k_n} $. The dimension calculation (\ref
  {dim-check}) above then shows that the term $ [\bar M_n (Y,\beta)]^{virt} $
  and hence the desired Gromov-Witten invariant will appear on the right hand
  side of our equation as one term among the $ D_k (X,A,B,M) $. The term
  coming from $ \bar M_\alpha (X,\beta) $ will vanish as $ \sum \alpha > d $,
  and all other terms are known recursively by exactly the same arguments as
  above for the relative invariants.
\end {proof}

\begin {remark}
  Although we have just shown that all restricted Gromov-Witten invariants of $
  Y \subset X $ can be computed from the Gromov-Witten invariants of $X$, only
  a very small subset of them is needed if one is only interested in the
  Gromov-Witten invariants of $Y$. First of all, analyzing the algorithm given
  in the proof above, one sees that it is sufficient to consider relative
  invariants of the form $ I_{(\alpha_1,0,\dots,0),\beta} (\gamma_1
  \psi_1^{k_1},\gamma_2, \dots,\gamma_n) $, i.e.\ we need multiplicities and
  cotangent line classes at only one of the marked points. In fact, in many
  cases it will be sufficient to look at invariants with only one marked point
  --- the WDVV equations of $Y$ can then be used to compute all Gromov-Witten
  invariants of $Y$. In a forthcoming paper we will give some explicit examples
  along these lines and show how corollary \ref {gw-hyp} can be used to reprove
  and generalize the ``mirror symmetry'' type formulas for Gromov-Witten
  invariants of certain hypersurfaces \cite {Be},\cite {G},\cite {LLY}.
\end {remark}



\begin {thebibliography}{XXXX}

\bibitem [B]{B} K. Behrend, \emph {Gromov-Witten invariants in algebraic
  geometry}, Inv.\ Math.\ \textbf {127} (1997), no.\ 3, 601--617.

\bibitem [Be]{Be} A. Bertram, \emph {Another way to enumerate rational curves
  with torus actions}, \preprint {math.AG}{9905159}.

\bibitem [BF]{BF} K. Behrend, B. Fantechi, \emph {The intrinsic normal cone},
  Inv.\ Math.\ \textbf {128} (1997), no.\ 1, 45--88.

\bibitem [BM]{BM} K. Behrend, Y. Manin, \emph {Stacks of stable maps and
  Gromov-Witten invariants}, Duke Math.\ J. \textbf {85} (1996), no.\ 1,
  1--60.

\bibitem [EGA4]{EGA4} A. Grothendieck, \emph {Éléments de géometrie
  algébrique IV}, IHES \textbf {32} (1967).

\bibitem [F]{F} W. Fulton, \emph {Intersection theory}, Springer 1984.

\bibitem [G]{G} A. Givental, \emph {Equivariant Gromov-Witten invariants},
  Internat.\ Math.\ Res.\ Notices \textbf {13} (1996), 613--663.

\bibitem [IP1]{IP1} E. Ionel, T. Parker, \emph {Gromov-Witten invariants of
  symplectic sums}, Math.\ Res.\ Lett.\ \textbf {5} (1998), no.\ 5, 563--576.

\bibitem [IP2]{IP2} E. Ionel, T. Parker, \emph {Relative Gromov-Witten
  invariants}, \preprint {math.SG}{9907155}.

\bibitem [K]{K} M. Kontsevich, \emph {Enumeration of rational curves via
  torus actions}, The moduli space of curves (Texel Island 1994),
  Birkhäuser Progr.\ in Math.\ \textbf {129} (1995), 335--368.

\bibitem [LLY]{LLY} B. Lian, K. Liu, S. Yau, \emph {Mirror principle I}, Asian
  J. of Math.\ \textbf {1} (1997), no.\ 4, 729--763.

\bibitem [LR]{LR} A. Li, Y. Ruan, \emph {Symplectic surgery and Gromov-Witten
  invariants of Calabi-Yau 3-folds I}, \preprint {math.AG}{9803036}.

\bibitem [P]{P} R. Pandharipande, \emph {Rational curves on hypersurfaces
  (after A. Givental)}, \preprint {math.AG}{9806133}.

\bibitem [R]{R} Y. Ruan, \emph {Surgery, quantum cohomology and birational
  geometry}, \preprint {math.AG}{9810039}.

\bibitem [S]{S} T. Shioda, \emph {Algebraic cycles on hypersurfaces in $ \PP^N
  $}, Algebraic geometry, Sendai 1985, Adv.\ Stud.\ Pure Math.\ \textbf {10}
  (1987), 717--732.

\bibitem [V]{V} R. Vakil, \emph {The enumerative geometry of rational and
  elliptic curves in projective space}, preprint available at
  \verb+http://www-math.mit.edu/~vakil/preprints.html+.

\end {thebibliography}
\end {document}